\documentclass{amsart}

\usepackage[utf8]{inputenc}
\usepackage{lmodern}
\usepackage[T1]{fontenc}
\usepackage{wasysym} 

\usepackage{amsmath,amsthm,amssymb}
\usepackage{mathtools}
\usepackage{amscd}

\usepackage{color}
\usepackage{hyperref}
\usepackage[inline]{enumitem}
\usepackage{array}
\usepackage{graphicx}
\usepackage{xcolor}
\usepackage{colonequals} 
\usepackage[normalem]{ulem}
\usepackage{mathrsfs}  

\usepackage[figuresleft]{rotating}
\usepackage{lscape}

\usepackage{blkarray}

\usepackage{romannum}
\AtBeginDocument{\pagenumbering{arabic}}

\usepackage{tkz-euclide,subfigure}

\usepackage{tikz}
\usepackage{tikz-cd}
\usepackage{tikz-3dplot}

\usetikzlibrary{shapes.geometric, calc}

\usetikzlibrary{decorations.markings,shapes.arrows, shapes.symbols}
    
\tikzset{paint/.style={ draw=#1!50!black, fill=#1!50 },     decorate with/.style={decorate,decoration={shape backgrounds,shape=#1,shape size=1mm}}}

\hypersetup{
    colorlinks=true, 
    linkcolor=blue, 
    urlcolor=red, 
    linktoc=all 
}

\newlist{enumarabic}{enumerate}{1}
\setlist[enumarabic]{font=\normalfont,label=(\arabic*),leftmargin=0.3in}
\newlist{enumroman}{enumerate}{1}
\setlist[enumroman]{font=\normalfont,label=(\roman*),leftmargin=0.3in}

\usepackage[all,cmtip]{xy}
\input xy
\xyoption{all}
\xyoption{rotate}
\xyoption{pdf}

\newcommand{\N}{\mathbb{N}}
\newcommand{\Q}{\mathbb{Q}}
\newcommand{\R}{\mathbb{R}}
\newcommand{\C}{\mathbb{C}}
\newcommand{\PP}{\mathbb{P}}
\newcommand{\Z}{\mathbb{Z}}

\newcommand{\Hom}{\mathrm{Hom}}
\newcommand{\xdownarrow}[1]{%
  {\left\downarrow\vbox to #1{}\right.\kern-\nulldelimiterspace}
}
\newcommand{\xuparrow}[1]{%
  {\left\uparrow\vbox to #1{}\right.\kern-\nulldelimiterspace}
}

\numberwithin{equation}{section} 

\DeclareMathOperator{\PD}{PD}

\makeatletter
\newtheorem*{rep@theorem}{\rep@title}
\newcommand{\newreptheorem}[2]{%
\newenvironment{rep#1}[1]{%
 \def\rep@title{#2 \ref{##1}}%
 \begin{rep@theorem}}%
 {\end{rep@theorem}}}
\makeatother 

\theoremstyle{plain}
\newtheorem{theorem}{Theorem}[section]
\newreptheorem{theorem}{Theorem}
\newtheorem{prop}[theorem]{Proposition}
\newtheorem{lemma}[theorem]{Lemma}

\newtheorem{cor}[theorem]{Corollary}

\theoremstyle{definition}
\newtheorem{dfn}[theorem]{Definition}
\newtheorem{eg}[theorem]{Example}

\newtheorem{remark}[theorem]{Remark}

\def\zk#1{\mathcal{Z}_{#1}}

\def\TMA{T^m}
\def\TMB{T^{m'}}
\def\TA{T^d}
\def\TB{T^{d'}}
\def\PA{P}
\def\PB{P'}
\def\lambdaA{\lambda}
\def\lambdaB{\lambda'}
\def\imathA{\imath}
\def\imathB{\imath'}

\def\sur{s}

\def\xq#1{X_{(#1)}}

\def\wSR{\mathsf{wSR}}
\def\SR{\mathsf{SR}}

\DeclareMathOperator{\Ima}{Im}

\DeclareMathOperator{\interior}{int}

\def\xabi{\xq{a_i,b_i}}

\def\mcup#1{M\big({#1}\big)}

\usepackage{chngcntr}

\makeatletter
\newcommand{\gettikzxy}[3]{%
  \tikz@scan@one@point\pgfutil@firstofone#1\relax
  \edef#2{\the\pgf@x}%
  \edef#3{\the\pgf@y}%
}
\makeatother

\makeatletter
\newcommand{\inlineitem}[1][]{%
\ifnum\enit@type=\tw@
    {\descriptionlabel{#1}}
  \hspace{\labelsep}%
\else
  \ifnum\enit@type=\z@
       \refstepcounter{\@listctr}\fi
    \quad\@itemlabel\hspace{\labelsep}%
\fi}
\makeatother

\newcolumntype{C}[1]{>{\centering\let\newline\\\arraybackslash\hspace{0pt}}m{#1}}

\begin{document}

\title[]
{The integral cohomology rings of four-dimensional toric orbifolds}
\author[X. Fu]{Xin Fu}
\address{Shanghai Institute for Mathematics and Interdisciplinary Sciences,  Shanghai 200433, China}
\email{x.fu@simis.cn}

\author[T. So]{Tseleung So}
\address{Institute of Mathematical Science, Pusan National University, Busan 46241, Republic of Korea}
\email{larry.so.tl@gmail.com}

\author[J. Song]{Jongbaek Song}
\address{Department of Mathematics Education, Pusan National University, Busan 46241, Republic of Korea}
\email{jongbaek.song@pusan.ac.kr}



\begin{abstract}%
Let $X(P,\lambda)$ be a 4-dimensional toric orbifold associated to a polygon $P$ and a characteristic function $\lambda$. Assuming that $X(P,\lambda)$ is locally smooth over a vertex of $P$, we determine the integral cohomology ring $H^*(X(P,\lambda);\Z)$ by constructing an explicit basis and expressing the cup products of the basis elements in terms of $P$ and $\lambda$.
\end{abstract}

\maketitle

\setcounter{tocdepth}{1}
\tableofcontents

\section{Introduction}
A toric orbifold is a $2d$-dimensional compact orbifold equipped with a locally standard $T^d$-action whose orbit space is a $d$-dimensional simple convex polytope~$P$. The $T^d$-action is encoded by a \emph{characteristic function} $\lambda \colon \mathcal{F}(P) \to \mathbb{Z}^d$ defined on the set~$\mathcal{F}(P)$ of facets of $P$ and satisfying certain conditions 
(See Definition~\ref{dfn_char_funct}).
In this setting, the toric orbifold is denoted by $X(P,\lambda)$.

A fundamental problem is to understand how the combinatorial data of $(P,\lambda)$ determine the topology of $X(P,\lambda)$, such as the ring structure of $H^*(X(P,\lambda))$.
In this direction, several cases have been studied. 
For example, when $X(P,\lambda)$ is smooth (that is, a \emph{quasitoric manifold}), there is a ring isomorphism
\[
H^*(X(P,\lambda);\Z)\cong \SR(P)/\mathcal{J}
\]
where $\SR(P)$ is the \emph{Stanley-Reisner ring} of $P$ and $\mathcal{J}$ is the ideal generated by the linear relations determined by $\lambda$ (See Lemma~\ref{lemma_J generator Phi_X bar}).
Without the smoothness assumption, this isomorphism still holds with rational coefficients. The reader is referred to~\cite{Dan, DJ, Jur} for more details. 

In general, very little is known about $H^*(X(P,\lambda);\Z)$ for a toric orbifold $X(P,\lambda)$. Nevertheless, assuming that $H^*(X(P,\lambda);\Z)$ is concentrated in even degrees, the work in~\cite{BSS, DKS} proves that there is a ring isomorphism
\begin{equation}\label{eqn_sing wSR/J}
H^*(X(P,\lambda);\Z)\cong \wSR(P,\lambda)/\mathcal{J},
\end{equation}
where $\wSR(P,\lambda)$ is a subring of $\SR(P)$ defined by the pair $(P,\lambda)$, called the \emph{weighted Stanley-Reisner ring}.
However, this description has limited practical use to determine cup products in $H^*(X(P,\lambda);\Z)$ 
in general because of the following reasons. 
First, it is difficult to write down explicit generators of $\wSR(P,\lambda)$. Second, 
even when all generators can be obtained from the definition, substantial effort is required to reduce them in order to extract an integral basis for $H^*(X(P,\lambda);\Z)$.

In this paper, we resolve the above problems for $4$-dimensional toric orbifolds by combining homotopy theoretic methods with the
algebraic isomorphism \eqref{eqn_sing wSR/J}. 
To describe our results, we first fix the notation. Let $P\subset\R^2$ be a polygon with~$n+2$ edges for $n\in\N$ and let $\lambda$ be a characteristic function on $P$. Label the vertices and edges of $P$ as described in Figure~\ref{fig_main char pair}, and write $\lambda(E_i)=(a_i,b_i)\in\Z^2$ for $1\leq i\leq n+2$.

\begin{figure}
\begin{tikzpicture}
\node[opacity=0, regular polygon, regular polygon sides=7, draw, minimum size = 2.2cm](m) at (0,0) {};
\coordinate (1) at (m.corner 1); \coordinate (2) at (m.corner 2); \coordinate (3) at (m.corner 3); 
\coordinate (4) at (m.corner 4); 
\coordinate (5) at (m.corner 5); 
\coordinate (6) at (m.corner 6); 
\coordinate (7) at (m.corner 7); 
\draw[fill=yellow!30, yellow!30] (1)--(2)--(3)--(4)--(5)--(6)--(7)--cycle;

\draw (1)--(2)--(3);
\draw (4)--(5);
\draw (6)--(7);

\draw (1)--($0.2*(7)+0.8*(1)$);
\draw (7)--($0.8*(7)+0.2*(1)$);

\draw (3)--($0.2*(4)+0.8*(3)$);
\draw (4)--($0.8*(4)+0.2*(3)$);

\draw (5)--($0.2*(6)+0.8*(5)$);
\draw (6)--($0.8*(6)+0.2*(5)$);

\foreach \a in {1,...,7} { 
\draw[fill] (\a) circle (1.5pt); 
}

\node[rotate=-28] at (m.side 7) {\scriptsize$\cdots$};
\node[rotate=130] at (m.side 3) {\scriptsize$\cdots$};
\node[rotate=51] at (m.side 5) {\scriptsize$\cdots$};

\node[left] at (m.side 2) {\scriptsize$E_{n+2}$};
\node[above left] at (m.side 1) {\scriptsize$E_{n+1}$};

\node[below] at (m.side 4) {\scriptsize$E_i$};
\node[right] at (m.side 6) {\scriptsize$E_j$};

\node[above] at (m.corner 1){\scriptsize$v_{n+1}$};
\node[left] at (m.corner 2){\scriptsize$v_{n+2}$};
\node[left] at (m.corner 3){\scriptsize$v_{1}$};
\node[below left] at (m.corner 4){\scriptsize$v_{i}$};
\node[below right] at (m.corner 5){\scriptsize$v_{i+1}$};
\node[below right] at (m.corner 6){\scriptsize$v_{j}$};
\node[right] at (m.corner 7){\scriptsize$v_{j+1}$};

\draw[->] (2,0)--(2.75,0);
\node[above] at (2.3,0) {$\lambda$};

\begin{scope}[xshift=150]
\node[opacity=0, regular polygon, regular polygon sides=7, draw, minimum size = 2.2cm](n) at (0.5,0) {};
\coordinate (1) at (n.corner 1); 
\coordinate (2) at (n.corner 2); 
\coordinate (3) at (n.corner 3); 
\coordinate (4) at (n.corner 4); 
\coordinate (5) at (n.corner 5); 
\coordinate (6) at (n.corner 6); 
\coordinate (7) at (n.corner 7); 
\draw[fill=yellow!30, yellow!30] (1)--(2)--(3)--(4)--(5)--(6)--(7)--cycle;

\draw (1)--(2)--(3);
\draw (4)--(5);
\draw (6)--(7);

\draw (1)--($0.2*(7)+0.8*(1)$);
\draw (7)--($0.8*(7)+0.2*(1)$);

\draw (3)--($0.2*(4)+0.8*(3)$);
\draw (4)--($0.8*(4)+0.2*(3)$);

\draw (5)--($0.2*(6)+0.8*(5)$);
\draw (6)--($0.8*(6)+0.2*(5)$);

\foreach \a in {1,...,7} { 
\draw[fill] (\a) circle (1.5pt); 
}
\end{scope}

\node[rotate=-28] at (n.side 7) {\scriptsize$\cdots$};
\node[rotate=130] at (n.side 3) {\scriptsize$\cdots$};
\node[rotate=51] at (n.side 5) {\scriptsize$\cdots$};

\node[left] at (n.side 2) {\scriptsize$(a_{n+2},b_{n+2})$};
\node[above left] at (n.side 1) {\scriptsize$(a_{n+1}, b_{n+1})$};

\node[below] at (n.side 4) {\scriptsize$(a_i,b_i)$};
\node[right] at (n.side 6) {\scriptsize$(a_j,b_j)$};
\end{tikzpicture}
\caption{A characteristic pair $(P,\lambda)$}
\label{fig_main char pair}
\end{figure}

We say that $v_{n+2}\in P$ is a \emph{smooth vertex} of $X(P,\lambda)$ if $\{\lambda(E_{n+1}), \lambda(E_{n+2})\}$ forms a basis for $\Z^2$ (see Definition~\ref{def_smooth vertex}).
Then by a suitable change of the basis of $\Z^2$, we assume
\begin{equation}\label{eqn_lambda refined form}
\lambda(E_{n+1})=(1,0)\quad
\text{and}\quad
\lambda(E_{n+2})=(0,1).
\end{equation}
Our main result determines the integral cohomology ring of such $X(P,\lambda)$ as follows.

\begin{theorem}\label{thm_main}
Let $X(P,\lambda)$ be a $4$-dimensional toric orbifold associated with $(P, \lambda)$ described as in Figure \ref{fig_main char pair}. If $v_{n+2}$ is a smooth vertex and $\lambda$ satisfies \eqref{eqn_lambda refined form},
then $\tilde{H}^*(X(P,\lambda))$ is torsion-free and has a basis
\[
\{u_1,\ldots,u_n;v\}
\]
with $\deg u_i=2$ for $i=1, \dots, n$ and $\deg v=4$ such that
\begin{equation}\label{eqn_cup prod equation}
u_i\cup u_j=a_ib_jv
\qquad
\text{for }1\leq i\leq j\leq n.
\end{equation}
\end{theorem}

The torsion-freeness of $\tilde{H}^*(X(P,\lambda))$ can be deduced from the work of~\cite{Fis,Jor,KMZ}. Our main contribution is 
an explicit construction of $u_i$'s and $v$ (see Definition \ref{dfn_good cellular basis}) together with a proof of Equation~\eqref{eqn_cup prod equation}, which will be given in Section \ref{sec_pf_main_thm}. 
To achieve this, we develop a theory of \emph{toric morphisms} for toric orbifolds of arbitrary dimensions. In particular, we study the ring homomorphisms of weighted Stanley--Reisner rings induced by toric morphisms. This novel approach is of independent interest as it applies to arbitrary toric orbifolds.

It is natural to compare Equation~\eqref{eqn_cup prod equation} with the intersection product in the case where $X(P,\lambda)$ is a toric surface.
Since cup products of degree two elements are commutative, both Equation~\eqref{eqn_cup prod equation} and the intersection product give rise to bilinear forms. In~\cite{FSS_chow ring}, we show that these two bilinear forms are inverse of one another, and explain how the generators $\{u_1,\ldots,u_n;v\}$ in Theorem~\ref{thm_main} correspond to divisor subvarieties. From this perspective, Equation~\eqref{eqn_cup prod equation} 
recovers the classical intersection theory in this special case, while also applying to a broader class of $4$-dimensional toric orbifolds.
In addition, it has motivated further  work on computing the signature and Steenrod operations for $4$-dimensional toric orbifolds. See~\cite{Mas} and~\cite{So}.

Theorem~\ref{thm_main} has several corollaries that compute cup products in $H^*(X(P,\lambda))$ for special cases in a more general setting. These results were included in Section~7 of the previous version of the paper. However, their proofs require additional assumptions and computation techniques that are different from those developed here. To keep this paper to a reasonable length and focus on the main ideas, we defer these results to a sequential paper~\cite{FSS_cup product two} currently in preparation.

This paper is organized as follows.
In Section~\ref{sec_preliminaries}, we review the combinatorial construction of toric orbifolds $X(P, \lambda)$ and moment-angle manifolds $\mathcal{Z}_P$, and relate their (equivariant) cohomology to weighted Stanley--Reisner rings~$\wSR(P,\lambda)$ and Stanley--Reisner rings $\SR(P)$. Several preliminary results on weighted Stanley--Reisner rings have also been proved.

In Section~\ref{sec_toric_morphism} we develop a general theory of  toric morphisms for toric orbifolds, motivated by toric morphisms of toric varieties~\cite[Chapter 3.3]{CLS} and of partial quotients~\cite[Section 3]{FF}. Section~\ref{sec_toric morphism in dim four} focuses on two special classes of toric morphisms for $4$-dimensional toric orbifolds: rescaling morphisms and edge-contraction morphisms.

In Section~\ref{sec_DTS}, we introduce degenerate toric spaces and define cellular bases for their cohomology. We further extend rescaling and edge-contraction morphisms to degenerate toric spaces, and prove several preparatory lemmas for later use.

In Section~\ref{sec_pf_main_thm} we prove Theorem~\ref{thm_main} by explicitly constructing a basis, called the \emph{algebraic cellular basis}, for $H^*(X(P,\lambda))$ and showing that the basis elements satisfy Equation~\eqref{eqn_cup prod equation}.

\medskip

\noindent 
\textbf{Glossary  of notation}. 
It may be helpful to list the main notational conventions used in the paper. 
\begin{itemize} 
\item $T^n$, an $n$-dimensional compact torus;
\item $P$, a $d$-dimensional simple polytope with $m$ facets (in Sections~\ref{sec_preliminaries} and~\ref{sec_toric_morphism}), or a polygon with $(n+2)$ edges (in Sections~\ref{sec_toric morphism in dim four} --~\ref{sec_pf_main_thm});
\item $\lambda, \Lambda$, a characteristic function and its representing matrix: Definition~\ref{dfn_char_funct} and Equation~\eqref{eq_char_matrix};
\item $T_E,S_E$, subtori of $T^d$ and $T^m$ determined by a face $E\subset P$: Equations~\eqref{eqn_T_E} and~\eqref{eqn_S_E};
\item $\mathcal{Z}_P$, moment-angle manifold: Definition~\ref{def_moment angle mfd};
\item $X(P,\lambda)$, a toric orbifold: Definition~\ref{dfn_toric_orb};
\item $X_{(\underline{a},\underline{b})}$, a $4$-dimensional toric orbifold or degenerate toric space with smooth vertex $v_{n+2}$: Equations~\eqref{eq_1001_char_fun} and~\eqref{eq_1001_deg_char_fun};
\item $\SR(P)$, the Stanley-Reisner ring of $P$: Definition~\ref{dfn_SR[P]}; 
\item $\wSR(P,\lambda)$, the weighted Stanley-Reisner ring of $(P,\lambda)$: Definition~\ref{dfn_wSR};
\item $\mathcal{J}$, the ideal generated by images of $\Ima(\varpi^*)$ in~\eqref{eqn_wSR bar isom prototype}: Lemma~\ref{lemma_J generator Phi_X bar};
\item $\overline{\wSR}(P,\lambda)$, the quotient of $\wSR(P,\lambda)$ by $\mathcal{J}$~\eqref{def_wSR bar};
\item $\Phi,\Phi_X,\overline{\Phi}_X$, isomorphisms from (equivariant) cohomology rings to (weighted) Stanley-Reisner rings and their quotients: see~\eqref{dfn_SR[P] isom prototype},~\eqref{eqn_wSR isom}, and~\eqref{eqn_wSR bar isom};
\item $(\psi_1,\psi_2),(\psi_1,\tilde{\psi}_2)$, a compatible pair and its lifting: Definitions~\ref{def_toric_morphism} and~\ref{dfn_lifting of compatible pair}
\item $X(\psi_1,\psi_2),\mathcal{Z}(\psi_1,\tilde{\psi}_2), \wSR(\psi_1,\psi_2), \overline{\wSR}(\psi_1,\psi_2),\SR(\psi_1,\tilde{\psi}_2)$, a toric morphism, its lifting, and induced ring homomorphisms: Definitions~\ref{def_toric_morphism} and~\ref{dfn_lifting of compatible pair}, and Equations~\eqref{eqn_wSR morphism}, \eqref{diagram_lifting and Z and X} and~\eqref{eqn_SR morphism lifting};
\item $g_{ij}$, the gcd function defined in~\eqref{eqn_rescaled characteristic function};
\item $X(id,\sigma_i)$, a rescaling morphism: Definitions~\ref{def_rescaling};
\item $X(\rho,id)$, an edge-contraction morphism: Definition~\ref{dfn_edge contracting morphism 1};
\item $\rho_i,\rho_{ij},\rho_{ij,k}$, edge-contractions from a polygon to a triangle, from a polygon to a square, and from a square to a triangle: Equations~\eqref{eqn_rho_i}, \eqref{eq_rho_ij}, and~\eqref{eqn_rho_ij,k};
\item $\{u_1,\ldots,u_n;v\},M(\mathcal{B},C_f)$, cellular basis and its cellular cup product representation: Definition~\ref{def_cellular_basis_cel_cup_repn}.
\end{itemize}

\subsection*{Acknowledgment}

Fu was supported by Basic Science Research Program through the National Research Foundation of Korea(NRF) funded by the Ministry of Education (NRF-2021R1A6A1A10044950 and NRF-2019R1A2C2010989), the National Natural Science Foundation of
China (no.\,12501083) and the startup fund of SIMIS. 
So was supported by Pacific Institute for the Mathematical Sciences (PIMS) Postdoctoral Fellowship, NSERC Discovery Grant and NSERC RGPIN-2020-06428.
So and Song are supported by the National Research Foundation of Korea(NRF) grant funded by the Korea government(MSIT) (RS-2025-00555914). Song was also supported by the KIAS Individual Grant (SP076101) at Korea Institute for Advanced Study.

The authors thank SIMIS for their hospitality and financial support for their research visit when this project was carried out.

\section{Preliminaries}\label{sec_preliminaries}
Let $S^1$ be the unit circle in $\mathbb{C}$ and $T^r=(S^1)^r$ the $r$-torus for any $r \in\mathbb{N}$. We identify $\mathbb{Z}^r$ with $\Hom(S^1, T^r)$, the lattice of $1$-parameter subgroups of $T^r$, and denote by $S^1_i$ the one-parameter subgroup of $T^r$ determined by the $i$-th unit vector in $\mathbb{Z}^r$.

\subsection{Toric orbifolds and moment-angle manifolds}\label{subsec_toricorb_MAC}
A $d$-dimensional polytope~$P$ is \emph{simple} if there are exactly $d$ facets (i.e., codimension-1 faces) intersecting at each vertex $v$. Let
\[
\mathcal{F}(P)=\{E_1,\ldots,E_m\}
\]
be the set of facets  in $P$.

\begin{dfn}\label{dfn_char_funct}
A \emph{characteristic pair} $(P,\lambda)$ is a $d$-dimensional simple polytope $P$ equipped with a map $\lambda\colon \mathcal{F}(P)\to\Z^d$,
called a \emph{characteristic function}, such that
\begin{enumerate}
\item
$\lambda(E_i)$ is a primitive vector for any $i=1,\dots, m$; 
\item
if $E_{i_1}\cap \cdots\cap E_{i_\ell} \neq \emptyset$, then $\{\lambda(E_{i_1}), \dots, \lambda(E_{i_\ell})\}$ is linearly independent.
\end{enumerate}
\end{dfn}

We call $\lambda(E_i)$ the \emph{characteristic vector} for $E_i$. Writing $\lambda(E_i)=(\lambda_{1i}, \dots, \lambda_{di})$, we often represent $\lambda$ by the following $(d\times m)$-matrix
\begin{equation}\label{eq_char_matrix}
\Lambda \colonequals \begin{pmatrix}
\lambda_{11}    &\lambda_{12}   &\cdots &\lambda_{1m}\\
\vdots          &\vdots         &\ddots &\vdots\\
\lambda_{d1}    &\lambda_{d2}   &\cdots &\lambda_{dm}
\end{pmatrix},
\end{equation}
called the \emph{characteristic matrix} of $\lambda$.

Let $E$ be a face of codimension-$\ell$ in $P$. Since $P$ is a simple polytope, we have
\[
E=E_{i_1}\cap \dots \cap E_{i_\ell}
\]
for some facets $E_{i_1},\ldots, E_{i_\ell}$ of $P$. In this case, we write
\begin{equation}\label{eqn_T_E}
T_E=\left\{\exp(t_1\lambda(E_{i_1})+\cdots+t_{\ell}\lambda(E_{i_{\ell}}))\mid t_1,\ldots,t_{\ell}\in\R\right\}
\end{equation}
as the rank-$\ell$ subtorus of $T^d$ generated by $\{\lambda(E_{i_1}), \dots, \lambda(E_{i_{\ell}})\}$
and write
\begin{equation}\label{eqn_S_E}
S_E=S^1_{i_1}\times\cdots\times S^1_{i_{\ell}}<\prod^m_{j=1}S^1_j.
\end{equation}
When $E=P$, we set $T_P$ as the trivial subgroup of $T^d$.

\begin{dfn}\label{dfn_toric_orb}
Given a characteristic pair $(P,\lambda)$, the associated \emph{toric orbifold} is the quotient space
\[
X(P, \lambda) \colonequals P \times T^d /_\sim
\]
equipped with a $T^d$-action and a map $\pi\colon X(P,\lambda)\to P$ defined as follows:
\begin{itemize}
\item
$(x,t)\sim (y,s)$ if and only if $x=y$ and $t^{-1}s\in T_{E(x)}$, where $E(x)$ is the minimal face of $P$ containing $x$ in its relative interior. We denote by  $[x,t]_{\sim}$ the equivalence class of $(x,t)\in P\times T^d$;
\item
the $T^d$-action $T^d\times X(P,\lambda)\to X(P,\lambda)$ is given by $(g,[x,t]_{\sim})\mapsto[x,gt]_{\sim}$ for any $g\in T^d$;
\item
the orbit map $\pi\colon X(P,\lambda)\to P$ is given by
$\pi([x,t]_{\sim})=x$.
\end{itemize}
\end{dfn}

\begin{remark}
A toric orbifold of Definition \ref{dfn_toric_orb} first appeared in \cite[Section~7]{DJ} and
Poddar--Sarkar gave an explicit axiomatic definition of a toric orbifold and proved that it is equivalent to the constructive definition given in Definition \ref{dfn_toric_orb}. We refer to \cite[Section 2]{PS}.
\end{remark}

An alternative construction of a toric orbifold $X(P, \lambda)$ is the quotient of the \emph{moment-angle manifold} $\mathcal{Z}_P$ defined below by a locally free action of a subtorus of~$T^m$ determined by $\lambda$. 

\begin{dfn}\label{def_moment angle mfd}
Let $P$ be a $d$-dimensional simple polytope with $m$ facets. Then the \emph{moment-angle manifold} associated with $P$ is the quotient space
\[
\mathcal{Z}_P\colonequals P\times T^m/_\approx,
\]
equipped with a $T^m$-action and a map $\tilde{\pi}\colon \mathcal{Z}_P\to P$ defined as follows:
\begin{itemize}
\item
$(x, t)\approx (y,s)$ if and only if $x=y$ and~$t^{-1}s \in S_{E(x)}$, and the equivalence class of $(x,t)\in P\times T^m$ is denoted by $[x,t]_{\approx}$;
\item
the $T^m$-action $T^m\times\mathcal{Z}_P\to\mathcal{Z}_P$ is given by
$(g,[x,t]_{\approx})\mapsto[x,g t]_{\approx}$ for any $g\in T^m$;
\item
the orbit map $\tilde{\pi}\colon \mathcal{Z}_P\to P$ is given by
$\tilde{\pi}([x,t]_{\approx})=x$.
\end{itemize}
\end{dfn}

The characteristic matrix $\Lambda$ of $\lambda$ in \eqref{eq_char_matrix} gives a short exact sequence 
\[
\begin{tikzcd}
    1 \rar &  \ker \exp \Lambda \rar &  T^m \rar{\exp \Lambda} & T^d \rar &  1
\end{tikzcd}
\]
where 
$(\exp\Lambda) (t_1, \dots, t_m) = \big(t_1^{\lambda_{11}}t_2^{\lambda_{12}}\cdots t_m^{\lambda_{1m}},\,\ldots, \,t_1^{\lambda_{d1}}t_2^{\lambda_{d2}}\cdots t_m^{\lambda_{dm}}\big)$. Then we have a~$T^d$-equivariant homeomorphism
\[
X(P, \lambda)\cong \mathcal{Z}_P/ \ker \exp \Lambda.
\]
We refer to \cite[Section 7]{DJ}. 
In what follows, we write
\begin{equation}\label{projection kappa}
\kappa \colon \mathcal{Z}_P \to X(P, \lambda)
\end{equation}
as the quotient map defined by $\kappa([x,t]_{\approx})=[x,(\exp\Lambda)(t)]_{\sim}$. 
Note that $\kappa$ is \emph{$\exp\Lambda$-equivariant}, namely, there is a commutative diagram
\[
\begin{tikzcd}
\mathcal{Z}_P\rar{\kappa}\dar{g}	&X(P,\lambda)\arrow[d,"(\exp\Lambda)(g)"]\\
\mathcal{Z}_P\rar{\kappa}				&X(P,\lambda)
\end{tikzcd}
\]
for any $g\in T^m$.

\subsection{Equivariant and ordinary cohomology rings}\label{subsec_equiv_coho}
For a Lie group~$G$ and a~$G$-space~$X$, we denote by $X_G\colonequals EG\times_G X$ the Borel construction of~$X$. As the map~$\kappa$ in \eqref{projection kappa} is $\exp \Lambda$-equivariant, the Borel fibration sequences of $\mathcal{Z}_P$ and $X(P,\lambda)$ fit into the following commutative diagram
\begin{equation}\label{diagram_borel fib}
\begin{tikzcd}
\mathcal{Z}_P\arrow{r}{\tilde{\imath}}\arrow{d}{\kappa}   &(\mathcal{Z}_P)_{T^m}\arrow{r}{\tilde{\varpi}}\arrow{d}{\kappa_T}   &BT^m\arrow{d}{B\exp\Lambda}\\
X(P,\lambda)\arrow{r}{\imath}  &X(P,\lambda)_{T^d}\arrow{r}{\varpi}   &BT^d,
\end{tikzcd}
\end{equation}
where $\kappa_T$ is induced by $\exp\Lambda$ and~$\kappa$.
The proof of Theorem~\ref{thm_main} is built on the fact that $H^*_T(\mathcal{Z}_P),H^*_{T^d}(X(P,\lambda))$ and~$H^*(X(P,\lambda))$ can be expressed in terms of degree~2 elements of $H^\ast (BT^m)$. 

\begin{dfn}\label{def_canonical_basis}
For any $r\in \mathbb{N}$, let $p_i\colon BT^r \to BS^1_i$ be the projection onto the $i$-th factor for $i=1, \dots, r$. Writing $x_i\colonequals p_i^\ast(u)\in H^\ast(BT^r)$ for the  standard generator~$u$ of $H^2(BS^1_i)$, we call $\{x_1, \dots, x_r\}$ the \emph{canonical basis} of $H^2(BT^r)$.
\end{dfn}

Given an $(r\times s)$-integral matrix
$A=(a_{ij})$, let $\exp A\colon T^s\to T^r$ be the induced Lie group homomorphism.
Let $\{x_1,\ldots,x_r\}$ and $\{y_1,\ldots,y_s\}$ be the canonical bases of~$H^2(BT^r)$ and $H^2(BT^s)$, respectively.  Then the homomorphism
\begin{equation}\label{eq_BexpA}
(B\exp A)^*\colon H^*(BT^r)\to H^*(BT^s)
\end{equation}
is given by
$(B\exp A)^*(x_i)=\sum^s_{j=1}a_{ij}y_j$.

We first study $H^*_{T^m}(\mathcal{Z}_P)$. Consider the Borel fibration sequence given in the first row of \eqref{diagram_borel fib}.
Let $\{y_1,\ldots,y_m\}$ be the canonical basis of $H^2(BT^m)$. Then $\widetilde{\varpi}$ induces a ring homomorphism
\begin{equation}\label{eq_tilde_varphi}
\widetilde{\varpi}^*\colon H^*(BT^m)=\Z[y_1,\ldots,y_m]\to H^*_{T^m}(\mathcal{Z}_P).
\end{equation}
It is known that $\widetilde{\varpi}^*$ is surjective and its kernel $\ker \widetilde{\varpi}^*$ is
\begin{equation*}
\mathcal{I} \colonequals \langle{y_{i_1}\cdots y_{i_k}\mid E_{i_1}\cap \cdots \cap E_{i_k} = \emptyset}\rangle,
\end{equation*}
see for instance \cite[Theorem 4.8]{DJ}.
Hence, there is an isomorphism of graded commutative rings
\begin{equation}\label{dfn_SR[P] isom prototype}
\Phi\colon H^*_{T^m}(\mathcal{Z}_P)\to \Z[y_1,\ldots,y_m]/\mathcal{I}.
\end{equation}
In the literature,
$\mathcal{I}$ is called the \emph{Stanley--Reisner ideal} and the quotient ring in~\eqref{dfn_SR[P] isom prototype} is called the \emph{Stanley--Reisner ring} of $P$, which we denote by $\SR(P)$. 

\begin{dfn}\label{dfn_SR[P]}
For the canonical basis $\{y_1,\ldots,y_m\}$ of $H^2(BT^m)$, we call the set 
\[
\{y_1+\mathcal{I},\ldots,y_m+\mathcal{I}\}\subset\SR(P)
\]
of quotient images the \emph{canonical generators} of $\SR(P)$. For simplicity, we keep using~$y_i$ for $y_i+\mathcal{I}\in \SR(P)$.
\end{dfn}

To study $H^\ast_T(X(P, \lambda))$, we begin with the following proposition inspired by the result of \cite[Lemma 5.6-(2)]{BSS}, where the authors used the language of fans of toric varieties.

\begin{prop}\label{subring}
Let $\kappa_T\colon (\mathcal{Z}_P)_{T^m}\to X(P,\lambda)_{T^d}$ be the map given in \eqref{diagram_borel fib}.
If the cohomology $H^*(X(P,\lambda))$ is concentrated in even degrees,
then the induced homomorphism 
\begin{equation}\label{eq_kappa_T}
\kappa_T^*\colon  H_{T^d}^\ast\big(X(P, \lambda)\big) \to H^*_{T^m}(\mathcal{Z}_P)
\end{equation}
is injective. 
\end{prop}

\begin{proof}
Let $V$ be the vertex set of $P$. For each $v\in V$, the restriction $\kappa_v\colonequals$  $\kappa|_{\tilde{\pi}^{-1}(v)}$ of $\kappa$ defined in \eqref{projection kappa} for the orbit map $\tilde{\pi}\colon \mathcal{Z}_P \to P$  maps~$\tilde{\pi}^{-1}(v)$ onto $\pi^{-1}(v)$. Hence we have a composition 
\[
\tilde{\pi}^{-1}(v)\overset{\kappa_v}{\longrightarrow}\pi^{-1}(v)\hookrightarrow X(P,\lambda),
\]
which gives us a commutative diagram
\begin{equation}\label{diagram_kappa restriction}
\begin{tikzcd}
\bigsqcup_{v\in V} \tilde{\pi}^{-1}(v) \arrow{d}[swap]{\bigsqcup_{v\in V} \kappa_v} \arrow[hook]{r} & \mathcal{Z}_P \arrow{d}{\kappa} \\
\bigsqcup_{v\in V}\pi^{-1}(v) \arrow[hook]{r} & X(P, \lambda).
\end{tikzcd}
\end{equation}
Notice that $\tilde{\pi}^{-1}(v)$ is a $T^m$-invariant subspace in $\mathcal{Z}_P$ and $\pi^{-1}(v)$ is a $T^d$-invariant subspace of $X(P, \lambda)$. Hence, we take the equivariant cohomology of~\eqref{diagram_kappa restriction}, 
which gives us 
\[
\begin{tikzcd}
\bigoplus_{v\in V}H^*_{T^m}\big(\tilde{\pi}^{-1}(v)\big) & H^*_{T^m}(\zk{P}) \lar \\
\bigoplus_{v\in V}H^*_{T^d}\big(\pi^{-1}(v)\big)
\uar{\bigoplus_{v\in V}(\kappa_v)_T^*}
&H^*_{T^d}\big(X(P,\lambda)\big). \lar[swap]{\gamma}\uar[swap]{\kappa_T^*}
\end{tikzcd}
\]
Since $H^\ast(X(P, \lambda))$ is concentrated in even degrees,  $\gamma$ is injective by \cite[Theorem~1.1]{FP}.
Hence, it suffices to show the injectivity of~$(\kappa_v)^*_T$ for all $v\in V$.

Fix $v\in V$ and consider the commutative diagram
\begin{equation}\label{diagram_kappa_v borel fib}
\begin{tikzcd}
\tilde{\pi}^{-1}(v)\rar \dar{\kappa_v} &\big(\tilde\pi^{-1}(v)\big)_{T^m}\rar{\tilde{\varpi}}\dar{(\kappa_v)_T} &BT^m\dar{B\exp\Lambda}\\
\pi^{-1}(v)\rar  &\big(\pi^{-1}(v)\big)_{T^d}\rar{\varpi}     &BT^d
\end{tikzcd}
\end{equation}
where the rows are Borel fibration sequences of $\tilde{\pi}^{-1}(v)$ and $\pi^{-1}(v)$. Since $P$ is a simple polytope of dimension $d$, there exist $d$ facets $E_{i_1},\dots,E_{i_d}$ intersecting $v$. 

Since $\tilde{\pi}^{-1}(v)\cong T^m/\prod^d_{j=1}S^1_{i_j}$ and $\pi^{-1}(v)\cong \{v\}$, we have 
\[
\big(\tilde\pi^{-1}(v)\big)_{T^m}\cong \prod^d_{j=1}BS^1_{i_j} \quad \text{and} \quad \big(\pi^{-1}(v)\big)_{T^d}\cong BT^d.
\]
Then $\tilde{\varpi}$ and $\varpi$ in~\eqref{diagram_kappa_v borel fib} can be identified with the inclusion $\prod^d_{j=1}BS^1_{i_j}\to BT^m$ and the identity map on $BT^d$, respectively. Thus, the right square of~\eqref{diagram_kappa_v borel fib} induces the commutative diagram
\begin{equation}\label{eq_cohom_vertex}
\begin{tikzcd}
H^*(BT^d)\rar{\varpi^\ast}\dar[swap]{(B\exp\Lambda)^*}  &H^*(BT^d)\dar{(\kappa_v)^*_T}\\
H^*(BT^m)\rar{\tilde{\varpi}^*} &H^*\big(\prod^d_{j=1}BS^1_{i_j}\big).
\end{tikzcd}
\end{equation}

Let $\{x_1,\ldots,x_d\}$ and $\{y_1,\ldots,y_m\}$ be canonical bases of $H^2(BT^d)$ and $H^2(BT^m)$, respectively. Then $\{y_{i_1}, \dots, y_{i_d}\}$ forms the canonical basis of $H^2(\prod_{j=1}^d BS_{i_j}^1)$. Moreover, writing $\lambda(E_i)=(\lambda_{1i},\ldots,\lambda_{di})$, we have
\[
\big(\tilde{\varpi}^*\circ(B\exp\Lambda)^*\big)(x_k)=\tilde{\varpi}^*\left(\sum^m_{i=1}\lambda_{ki}\,y_{i}\right)=\sum^d_{j=1}\lambda_{k,i_j}y_{i_j},
\]
where the first equality follows from \eqref{eq_BexpA} and  the second one follows because $\tilde{\varpi}^\ast$ maps $y_i$ to itself for $i\in \{i_1, \dots, i_d\}$ and zero otherwise. Hence, the commutativity of \eqref{eq_cohom_vertex} implies that the matrix representation of $(\kappa_v)^*_T$ is the following square matrix 
\[
\begin{pmatrix}
\lambda(E_{i_1})^t & \cdots & \lambda(E_{i_d})^t
\end{pmatrix}.
\]
As $\{\lambda(E_{i_1}),\ldots,\lambda(E_{i_d})\}$ is linearly independent due to Definition~\ref{dfn_char_funct}, the map  
$(\kappa_v)^*_T$ is injective for each $v$. 
Hence, we conclude that $\kappa^*_T$ is injective as well.
\end{proof}

\begin{dfn}\label{dfn_wSR}
We say a characteristic pair $(P,\lambda)$ is \emph{even} if $H^*(X(P,\lambda))$ is concentrated in even degrees.
For such $(P,\lambda)$, its \emph{weighted Stanley--Reisner ring} is the subring
\[
\wSR(P,\lambda)\colonequals \Ima(\Phi\circ \kappa^*_T)
\]
of $\SR(P)$, where $\Phi$ and $\kappa_T^\ast$ are defined in \eqref{dfn_SR[P] isom prototype} and \eqref{eq_kappa_T}, respectively.
\end{dfn}

\begin{remark}
The notion of a weighted Stanley--Reisner ring first appeared in~\cite{BSS} for projective toric orbifolds corresponding to \emph{even} characteristic pairs.  
Darby--Kuroki--Song generalized it to a much broader class of orbifolds known as torus orbifolds, where they call it the \emph{weighted face ring}~\cite[Definition 3.2]{DKS} of the GKM-graph associated with the given torus orbifold. 
\end{remark}

For an even characteristic pair $(P,\lambda)$, 
we denote by 
\begin{equation}\label{eqn_wSR isom}
\Phi_X\colon H^*_{T^d}(X(P,\lambda))\to\wSR(P,\lambda)
\end{equation}
the isomorphism $\Phi\circ \kappa^*_T$ onto its image. 
Proposition~\ref{subring} and Definition~\ref{dfn_wSR} imply that the following diagram commutes
\begin{equation}\label{diagram_Phi Phi_X commute}
\begin{tikzcd}
H^*_{T^d}(X(P,\lambda))\rar{\Phi_X}\dar[swap]{\kappa^*_T} &\wSR(P,\lambda)\dar[hook]\\
H^*_T(\mathcal{Z}_P)\rar{\Phi} &\SR(P).
\end{tikzcd}
\end{equation}

We now study $H^*(X(P,\lambda))$ using \eqref{diagram_Phi Phi_X commute}. Consider the Borel fibration  
given in the second row of \eqref{diagram_borel fib}, which induces a sequence 
\[
H^*(BT^d)\overset{\varpi^*}{\longrightarrow}H^*_{T^d}(X(P,\lambda))\overset{\imath^*}{\longrightarrow}H^*(X(P,\lambda)).
\]
Since $H^\ast(X(P, \lambda))$ is concentrated in even degrees, the Eilenberg--Moore spectral sequence collapses at $E_2$-page and 
$\imath^*$ induces a ring isomorphism
\begin{equation}\label{eqn_wSR bar isom prototype}
H^*_{T^d}(X(P,\lambda))/\langle\Ima(\varpi^*)\rangle\to H^*(X(P,\lambda))
\end{equation}
where $\langle\Ima(\varpi^*)\rangle$ is the ideal generated by the image of $\varpi^*$. For the isomorphism~$\Phi_X$ of \eqref{eqn_wSR isom}, let $\mathcal{J}\colonequals \Phi_X\big(\langle\text{Im}(\varpi^*)\rangle\big)$ be the corresponding ideal in
$\wSR(P,\lambda)$ and write
\begin{equation}\label{def_wSR bar}
\overline{\wSR}(P,\lambda)\colonequals\wSR(P,\lambda)/\mathcal{J}.
\end{equation}
Then $\Phi_X$ descends to the following graded ring isomorphism 
\begin{equation}\label{eqn_wSR bar isom}
\overline{\Phi}_X\colon H^*(X(P,\lambda))\to \overline{\wSR}(P,\lambda).
\end{equation}

\begin{remark}\label{rmk_toric_manifold}
For a toric manifold $X(P,\lambda)$, the map $\kappa_T$ of Proposition \ref{subring} is a homotopy equivalence (see~\cite[Section 4.1]{DJ}). Hence, it follows that 
\[
\wSR(P,\lambda)=\SR(P)
\quad\text{and}\quad
\overline{\wSR}(P,\lambda)=\SR(P)/\mathcal{J}.
\]
\end{remark}

In the following lemma, we study the generators of the ideal $\mathcal{J}$ more explicitly.

\begin{lemma}\label{lemma_J generator Phi_X bar}
Let $(P,\lambda)$ be an even characteristic pair with its associated characteristic matrix defined in~\eqref{eq_char_matrix}.  Let $\{y_1,\ldots,y_m\}$ be the canonical generators of~$\SR(P)$. Then the ideal $\mathcal{J}$ is generated by the following linear terms
\[
\{\lambda_{k1}y_1 + \cdots + \lambda_{km}y_m
\mid k=1,\ldots,d\}.
\]
Furthermore, there is a commutative diagram
\begin{equation}\label{diagram_Phi_X Phi_X bar commute}
\begin{tikzcd}
H^*_{T}(X(P, \lambda))\rar{\imath^*}\dar[swap]{\Phi_X}  & H^*(X(P,\lambda))\dar{\overline{\Phi}_X}\\
\wSR(P,\lambda)\rar[two heads]            &\overline{\wSR}(P,\lambda)
\end{tikzcd}
\end{equation}
where the bottom arrow is the quotient map.
\end{lemma}

\begin{proof}
The proof of the first claim is similar to the one given in \cite[Proof of Theorem 5.3]{BSS}. 
Here we provide a proof in the context of this paper for readers' convenience. 
Indeed, it suffices to show that
\[
(\Phi_X\circ\varpi^*)(x_k)=\lambda_{k1}y_1+\cdots+\lambda_{km}y_m
\]
for $k=1,\ldots,d$, where $\{x_1,\ldots,x_d\}$ and $\{y_1, \dots, y_m\}$ are the canonical bases of~$H^2(BT^d)$ and $H^2(BT^m)$, respectively. 
Consider the commutative diagram
\[
\begin{tikzcd}
H^*(BT^d)\arrow{rr}{(B\exp\Lambda)^*}\dar{\varpi^*}   &&H^*(BT^m)\rar{\cong}\dar{\tilde{\varpi}^*}   &\Z[y_1, \dots, y_m]\dar[two heads]{q}\\
H_{T^d}^\ast(X(P, \lambda))\arrow{rr}{\kappa^*_T} && H^*_{T^m}(\mathcal{Z}_P)\rar{\Phi}   &\SR(P)
\end{tikzcd}
\]
where the left square is induced from the right square of \eqref{diagram_borel fib} and $q$ is the quotient map.
We note that the right square commutes due to the definition of $\Phi$ in~\eqref{dfn_SR[P] isom prototype}. 
As $\Phi_X=\Phi\circ\kappa^*_T$, we have
\[
 (\Phi_X\circ\varpi^\ast)(x_k)
=\left(q\circ(B\exp\Lambda)^\ast \right)(x_k)=\lambda_{k1}y_1+\cdots+\lambda_{km}y_m,   
\]
which establishes the first claim. 

Next we prove the commutativity of~\eqref{diagram_Phi_X Phi_X bar commute}. Consider the diagram
\[
\begin{tikzcd}
H^*_{T^d}(X(P,\lambda))\rar[two heads]\dar{\Phi_X} &H^*_{T^d}(X(P,\lambda))/\left<\Ima(\varpi^*)\right> \rar{\cong}\dar{f} &H^*(X(P,\lambda))\dar{\overline{\Phi}_X}\\
\wSR(P,\lambda)\rar[two heads] &\wSR(P,\lambda)/\mathcal{J}\rar{=}    &\overline{\wSR}(P,\lambda)
\end{tikzcd}
\]
where the horizontal maps in the left square are quotient maps and $f$ is induced by $\Phi_X$. The left square commutes due to the definition of $\mathcal{J}$ and the right square commutes due to the definition of $\overline{\Phi}_X$. Hence, the outer rectangle commutes. Since the composite of maps in the top row is $\imath^*\colon H^*_{T^d}(X(P,\lambda))\to H^*(X(P,\lambda))$, we obtain~\eqref{diagram_Phi_X Phi_X bar commute}.
\end{proof}

\begin{eg}\label{eg_Xab}
Consider the characteristic pair $(P, \lambda_{(a,b)})$ where $P$ is a triangle~$\triangle$ and $\lambda_{(a,b)}$ is the characteristic function on $P$ in Figure \ref{fig_char_pair_tri} for some nonzero integers~$a$ and $b$. 
\begin{figure}
\begin{tikzpicture}[ baseline=(current bounding box.center)]
\node[regular polygon, regular polygon sides=3, draw, minimum size = 1.5cm](Q) at (0,0) {};
\node[below] at (Q.side 2) {\scriptsize$E_1$}; 
\node[left] at (Q.side 1) {\scriptsize$E_3$}; 
\node[right] at (Q.side 3) {\scriptsize$E_2$}; 

\draw[->] (1.8,0)--(2.8,0);
\node[above] at (2.3,0) {$\lambda_{(a,b)}$};

\node[regular polygon,  regular polygon sides=3, draw, minimum size = 1.5cm](Qi) at (5,0) {};
\node[below] at (Qi.side 2) {\scriptsize$(a,b)$}; 
\node[left] at (Qi.side 1) {\scriptsize$(0,1)$}; 
\node[right] at (Qi.side 3) {\scriptsize$(1,0)$}; 
\end{tikzpicture}
\caption{Characteristic function on $\triangle$.}
\label{fig_char_pair_tri}
\end{figure}
The characteristic matrix associated $\lambda_{(a,b)}$ is given by 
\[
\left(\begin{array}{c c c}
a   &1  &0\\
b   &0  &1
\end{array}\right).
\]
Let $\{y_1,y_2,y_3\}$ be the canonical generators of $\SR(\triangle)$.  
Then, Lemma~\ref{lemma_J generator Phi_X bar} implies that the ideal $\mathcal{J}\subset \wSR(\triangle,\lambda_{(a,b)})$ is generated by~$\{ay_1+y_2,by_1+y_3\}$. 
The computational result of \cite[Example 5.4]{BSS} shows that $\overline{\wSR}(\triangle,\lambda_{(a,b)})$ is spanned by~$\{ 1,[aby_1],[y_2y_3]\}$ with the product structure $[aby_1]^2=ab[y_2y_3]$, where $[aby_1]$ and $[y_2y_3]$ are the quotient images of $aby_1$ and $y_2y_3$ modulo $\mathcal J$.
\end{eg}

\section{Toric morphisms}\label{sec_toric_morphism}
A toric morphism between toric varieties refers to an equivariant map with respect to the associated torus actions. In this section, we define a toric morphism in the context of toric orbifolds discussed in Section \ref{sec_preliminaries} and develop tools to calculate integral cohomology rings of toric orbifolds.

\begin{dfn}\label{def_toric_morphism}
Let $(P,\lambda)$ and $(P',\lambda')$  be characteristic pairs of dimension $d$ and~$d'$, respectively. Then
\begin{itemize}
\item
a \emph{compatible pair} $(\psi_1,\psi_2)\colon \PA\times \TA\to \PB\times \TB$ is the product map of a continuous map $\psi_1\colon P \to P'$ and a Lie group homomorphism $\psi_2 \colon T^d \to T^{d'}$ such that for any face $E$ of $P$ there exists a face $E'$ of $P'$  satisfying
\begin{equation}\label{eqn_compatible pair_E, E'}
\psi_1(E)\subseteq E'\quad
\text{and}\quad
\psi_2(T_{E}) \leq T_{E'}
\end{equation}
where $T_E$ and $T_{E'}$ are subtori given in~\eqref{eqn_T_E};
\item
the \emph{toric morphism} associated with a compatible $(\psi_1,\psi_2)$ is the map
\[
X(\psi_1,\psi_2)\colon X(\PA,\lambdaA)\to X(\PB,\lambdaB)
\]
given by $X(\psi_1,\psi_2)([x,t]_{\sim})=[\psi_1(x),\psi_2(t)]_{\sim}$ for $x\in \PA$ and $t\in \TA$.
\end{itemize}
\end{dfn}

Notice that $X(\psi_1,\psi_2)$ is $\psi_2$-equivariant, so there is a commutative diagram
\[
\begin{tikzcd}
X(\PA,\lambdaA)\rar{\imathA}\dar{X(\psi_1,\psi_2)}  &X(\PA,\lambdaA)_{\TA}\rar{\varpi}\dar{X(\psi_1,\psi_2)_T}  &B\TA\dar{B\psi_2}\\
X(\PB,\lambdaB)\rar{\imathB} &X(\PB,\lambdaB)_{\TB}\rar{\varpi'}  &B\TB
\end{tikzcd}
\]
of Borel fibrations,
where $X(\psi_1,\psi_2)_T([e,z])=[E\psi_2(e),X(\psi_1,\psi_2)(z)]$ for~$e\in E\TA$ and $z\in X(\PA,\lambdaA)$.
The left square gives a commutative diagram
\begin{equation}\label{diagram_X(psi_1 psi_2)^* n X(psi_1 psi_2)^*_T}
\begin{tikzcd}
H^*_{\TB}(X(\PB,\lambdaB))\arrow{rr}{X(\psi_1,\psi_2)^*_T}\arrow{d}[left]{(\imathB)^*} &&H^*_{\TA}(X(\PA,\lambdaA))\dar{\imathA^*}\\
H^*(X(\PB,\lambdaB))\arrow{rr}{X(\psi_1,\psi_2)^*} &&H^*(X(\PA,\lambdaA)).
\end{tikzcd}
\end{equation}

Suppose $(P,\lambda)$ and $(P',\lambda')$ are even (see Definition \ref{dfn_wSR}). Then the ring isomorphisms $\Phi_X$ and $\overline{\Phi}_X$ studied in~\eqref{eqn_wSR isom} and~\eqref{eqn_wSR bar isom} give us the following composites of ring homomorphisms
\begin{align}\label{eqn_wSR morphism}
\begin{split}
\wSR(\psi_1,\psi_2)&\colonequals\Phi_{X}\circ X(\psi_1,\psi_2)^*_T\circ\Phi_{X'}^{-1}\colon \wSR(\PB,\lambdaB)\to\wSR(\PA,\lambdaA), \\ 
\smallskip
\overline{\wSR}(\psi_1,\psi_2)&\colonequals\overline{\Phi}_{X}\circ X(\psi_1,\psi_2)^*\circ\overline{\Phi}_{X'}^{-1}\colon \overline{\wSR}(\PB,\lambdaB)\to\overline{\wSR}(\PA,\lambdaA).
\end{split}
\end{align}
Then, Diagram \eqref{diagram_X(psi_1 psi_2)^* n X(psi_1 psi_2)^*_T} can be extended to the following cubical diagram
\begin{equation}\label{diagram_naturality of wSR n wSR bar}
\small
\begin{tikzcd}
    &H^*_{\TA}(X(\PA,\lambdaA))
\arrow{rr}{\Phi_{X}}
\ar[]{dd}[swap,near end]{\imathA^*}
  &   &
\wSR(\PA,\lambdaA)
\arrow[dd, two heads]
\\
H^*_{\TB}(X(\PB,\lambdaB))
\ar[crossing over]{rr}[near end]{\Phi_{X'}}
\ar{dd}[swap]{(\imathB)^*}
\ar{ur}{X(\psi_1,\psi_2)^*_T}   &
&
\wSR(\PB,\lambdaB)
\ar{ur}[swap,near start]{\wSR(\psi_1,\psi_2)}  &\\
    &H^*(X(\PA,\lambdaA))
\ar{rr}[near start]{\overline{\Phi}_{X}}
  &\rar{}   &\overline{\wSR}(\PA,\lambdaA)\\
H^*(X(\PB,\lambdaB))
\ar{rr}{\overline{\Phi}_{X'}}
\ar{ur}{X(\psi_1,\psi_2)^*}   &    &
\overline{\wSR}(\PB,\lambdaB)
\ar[crossing over,twoheadleftarrow]{uu}
\ar{ur}[swap,near start]{\overline{\wSR}(\psi_1,\psi_2)},  &
\end{tikzcd}
\end{equation}
such that all faces commute, where the unnamed vertical arrows are quotient maps. Indeed, the left face is~\eqref{diagram_X(psi_1 psi_2)^* n X(psi_1 psi_2)^*_T}. The top and the bottom faces commute due to the definitions of $\wSR(\psi_1,\psi_2)$ and $\overline{\wSR}(\psi_1,\psi_2)$. The front and the rear faces commute due to~\eqref{diagram_Phi Phi_X commute} and \eqref{diagram_Phi_X Phi_X bar commute}. Since $\Phi_X$ and $\overline{\Phi}_X$ are isomorphisms, a diagram chasing shows that the right face of the cube also commutes.

The commutativity of \eqref{diagram_naturality of wSR n wSR bar} leads us to study the induced maps $X(\psi_1,\psi_2)^*$ and~$X(\psi_1,\psi_2)^*_T$ via $\wSR(\psi_1,\psi_2)$ and $\overline{\wSR}(\psi_1,\psi_2)$. The rest of this section is devoted to derive explicit formulas to calculate $\wSR(\psi_1,\psi_2)$ and $\overline{\wSR}(\psi_1,\psi_2)$, which will be given in Theorem \ref{thm_SR_morphisms_formula} and Corollary \ref{cor_wSR_n_wSR_bar_morphisms_formula}.

\begin{dfn}\label{dfn_lifting of compatible pair}
Let $(P,\lambda)$ and $(P',\lambda')$ be two characteristic pairs as in Definition~\ref{def_toric_morphism}.  Let $m$ and $m'$ be the numbers of facets in $P$ and $P'$, respectively. A \emph{lifting} of a compatible pair $(\psi_1,\psi_2)$ is the product map
\[
(\psi_1,\tilde{\psi_2})\colon \PA\times \TMA\to \PB\times \TMB
\]
 of $\psi_1$ and a Lie group homomorphism $\tilde{\psi}_2\colon \TMA \to \TMB$ 
satisfying the following conditions:
\begin{enumerate}
\item\label{dfn_lifting of compatible pair subgp condition}
if faces $E\subset P$ and $E'\subset P'$ satisfies~\eqref{eqn_compatible pair_E, E'} then
\[
\tilde{\psi_2}(S_E)\leq S_{E'}
\]
where $S_E$ and $S_{E'}$ are subtori given in~\eqref{eqn_S_E};

\item\label{dfn_lifting of compatible pair comm diagram}
There is a commutative diagram
\begin{equation}\label{eq_lift_comm}
\begin{tikzcd}
\TMA\rar{\tilde{\psi}_2}\arrow{d}[left]{\exp\Lambda}  &\TMB\dar{\exp\Lambda'}\\
\TA\rar{\psi_2}        &\TB
\end{tikzcd}
\end{equation}
\end{enumerate}
where $\Lambda$ and $\Lambda'$ are the characteristic matrices of $\lambda$ and $\lambda'$, respectively. 
\end{dfn}

Proposition \ref{prop_uniqueness_lifting} below shows that the lifting of a compatible pair $(\psi_1,\psi_2)$ is unique if it exists. However, not every compatible pair admits a lifting. See Example~\ref{ex_nonexistence_lifting}.

\begin{prop}\label{prop_uniqueness_lifting}
Let $(\psi_1,\psi_2)\colon (\PA,\lambdaA) \to (\PB,\lambdaB)$ be a compatible pair.
If a lifting of $(\psi_1,\psi_2)$ exists, then it is unique.
\end{prop}

\begin{proof}
Let $(\psi_1, \tilde{\psi}_2)$ be a lifting of $(\psi_1, \psi_2)$. 
For each facet $E_j\in \mathcal{F}(P)$, since $(\psi_1, \psi_2)$ is a compatible pair (see Definition \ref{def_toric_morphism}), there is a face $E'$ of $P'$
such that
\[\text{$\psi_1(E_j)\subset E'$ and $\psi_2(T_{E_j})\leq T_{E'}$}.\]
Note that $T_{E_j}$ is a one-parameter subgroup of~$\TA$ generated by~$\lambdaA(E_j)$, so $\psi_2(T_{E_j})$ is a subtorus in $T_{E'}$ of rank $\leq 1$.
Writing $\Psi$ for the $(d'\times d)$-integral matrix such that $\psi_2=\exp \Psi$, we have
\[
\Psi(\lambdaA(E_j))= c_{j_1}\lambdaB(E'_{j_1}) + \cdots + c_{j_\ell} \lambdaB(E'_{j_\ell}),
\]
where $\{j_1,\ldots,j_\ell\}$ is the unique set of indices such that $E'=E'_{j_1} \cap \cdots \cap E'_{j_\ell}$ is the intersection of facets of $P'$.  
Since $\{ \lambdaB(E'_{j_1}), \dots, \lambdaB(E'_{j_\ell})\}$ is linearly independent (see Definition \ref{dfn_char_funct}), the set $\{c_{j_1}, \dots, c_{j_k}\}$ is unique.

Since $(\psi_1, \tilde{\psi}_2)$ is a lifting of $(\psi_1, \psi_2)$, following Definition~\ref{dfn_lifting of compatible pair} we have 
\begin{equation}\label{eq inclusion}
 \tilde{\psi}_2(S_{E_j})\leq S_{E'}. 
\end{equation}
Let $\tilde{\Psi}=(d_{ij})$ be the $(m'\times m)$-integral matrix such that $\tilde{\psi}_2=\exp \tilde{\Psi}$. For each~$j$, we  have $d_{ij}=0$ if $i\notin\{j_1,\ldots, j_\ell\}$ due to the above inclusion~\eqref{eq inclusion}. If $i\in\{j_1,\ldots, j_\ell\}$, the commutativity of \eqref{eq_lift_comm} shows
$d_{ij}=c_{i}$. Finally, the result follows by the uniqueness of $c_{j_k}$'s.
\end{proof}

Note that the lifting $(\psi_1,\tilde{\psi}_2)$ defines a map
$\mathcal{Z}(\psi_1,\tilde{\psi}_2)\colon \mathcal{Z}_{\PA}\to\mathcal{Z}_{\PB}$
by
\[
\mathcal{Z}(\psi_1,\tilde{\psi}_2)([x,t]_{\approx})=[\psi_1(x),\tilde{\psi}_2(t)]_{\approx}
\]
for $x\in \PA$ and $t\in \TMA$.
By Definition~\ref{dfn_lifting of compatible pair}, it makes the following diagram commute
\begin{equation}\label{diagram_lifting and Z and X}
\begin{tikzcd}
\mathcal{Z}_{\PA}\arrow{rr}{\mathcal{Z}(\psi_1,\tilde{\psi}_2)}\arrow{d}[left]{\kappa}   
&&\mathcal{Z}_{\PB}\dar{\kappa'}\\
X(\PA,\lambdaA)\arrow{rr}{X(\psi_1,\psi_2)}    &&X(\PB,\lambdaB)
\end{tikzcd}
\end{equation}
where $\kappa$ and $\kappa'$ are the quotient maps discussed in~\eqref{dfn_char_funct}. As the map $\mathcal{Z}(\psi_1,\tilde{\psi}_2)$ is~$\tilde{\psi}_2$-equivariant, the above \eqref{diagram_lifting and Z and X} can be extended to the following cubical diagram
\begin{equation}\label{diagram_naturality SR map proof cube}
\small
\begin{tikzcd}
    &B\TMA
\arrow{rr}{B\tilde{\psi}_2} 
\ar{dd}[swap,near end]{B\Lambda}
  &   &
B\TMB
\ar{dd}{B\Lambda'}
\\
(\mathcal{Z}_{\PA})_{\TMA}
\ar[crossing over]{rr}[near end]{\mathcal{Z}(\psi_1,\tilde{\psi}_2)_T}
\ar{ur}{\tilde{\varpi}}
\ar{dd}[swap]{\kappa_T}   &
&(\mathcal{Z}_{\PB})_{\TMB}
\ar{ur}[swap]{\tilde{\varpi}'}
&\\
&B\TA
\ar{rr}[near start]{B\psi_2}   
&   &B\TB,\\
X(\PA,\lambdaA)_{\TA}
\ar{rr}{X(\psi_1,\psi_2)_T}\ar{ur}{\varpi} 
&    &
X(\PB,\lambdaB)_{\TB}\ar{ur}[swap]{\varpi'} 
\ar[crossing over,leftarrow]{uu}[swap,near end]{\kappa'_T}&
\end{tikzcd}
\end{equation}
where $\mathcal{Z}(\psi_1,\tilde{\psi}_2)_T$ is the map induced by  $E\tilde{\psi}_2$ and $\mathcal{Z}(\psi_1,\tilde{\psi}_2)$. Applying the cohomology to the front face of~\eqref{diagram_naturality SR map proof cube}, we have a commutative diagram
\begin{equation}\label{diagram_SR(psi_1 psi_2) definition}
\begin{tikzcd}
H^*_{\TB}(X(\PB,\lambdaB))\arrow{rr}{X(\psi_1,\psi_2)^*_T}\arrow{d}[left]{(\kappa'_T)^*} &&H^*_{\TA}(X(\PA,\lambdaA))\dar{(\kappa_T)^*}\\
H^*_{\TMB}(\mathcal{Z}_{\PB})\arrow{rr}{\mathcal{Z}(\psi_1,\tilde{\psi}_2)^*_T} &&H^*_{\TMA}(\mathcal{Z}_{\PA}).
\end{tikzcd}
\end{equation}

Now we define a ring homomorphism
\begin{equation}\label{eqn_SR morphism lifting}
\SR(\psi_1,\tilde{\psi}_2)\colonequals\Phi\circ\mathcal{Z}(\psi_1,\tilde{\psi}_2)^*_T\circ(\Phi')^{-1}\colon \SR(\PB)\to\SR(\PA)
\end{equation}
using isomorphisms $\Phi\colon H^*_{\TMA}(\mathcal{Z}_P)\to\SR(P)$ and $\Phi'\colon H^*_{\TMB}(\mathcal{Z}_{P'})\to\SR(P')$ studied in~\eqref{dfn_SR[P] isom prototype}. Then~\eqref{diagram_SR(psi_1 psi_2) definition} can be extended to the following cubical diagram
\begin{equation}\label{diagram_naturality SR map cube}
\small
\begin{tikzcd}
   &H^*_{\TA}(X(\PA,\lambdaA))
    \ar{rr}{\Phi_{X}} 
    \ar[]{dd}[swap,near end]{\kappa^*_T}
 &   &\wSR(\PA,\lambdaA)
\\
H^*_{\TB}(X(\PB,\lambdaB))
\ar[crossing over,near end]{rr}{\Phi_{X'}}
\ar{ur}{X(\psi_1,\psi_2)^*_T}
\ar{dd}[swap]{(\kappa'_T)^*}  
&    &\wSR(\PB,\lambdaB)
\ar{ur}[swap,near start]{\wSR(\psi_1,\psi_2)}
  &\\
    &H^*_{\TMA}(\mathcal{Z}_{\PA})
\ar[near start]{rr}{\Phi}  
&  &\SR(\PA)\ar[hookleftarrow]{uu}
\\
H^*_{\TMB}(\mathcal{Z}_{\PB})
\ar{rr}{\Phi'}\ar{ur}{\mathcal{Z}(\psi_1,\tilde{\psi}_2)^*_T}   &    &\SR(\PB)\ar{ur}[swap,near start]{\SR(\psi_1,\tilde{\psi}_2)}
\ar[hookleftarrow, crossing over]{uu}&
\end{tikzcd}
\end{equation}
such that all faces commute, where the unnamed vertical maps are inclusions.

In the following theorem, we show how the morphism $\mathcal{Z}(\psi_1,\tilde{\psi}_2)$ induces $\SR(\psi_1,\tilde{\psi}_2)$, which in turn determines $\wSR(\psi_1,\psi_2)$ and $\overline{\wSR}(\psi_1,\psi_2)$ in Corollary~\ref{cor_wSR_n_wSR_bar_morphisms_formula}.

\begin{theorem}\label{thm_SR_morphisms_formula}
Let $(P,\lambda)$ and $(P',\lambda')$ be even characteristic pairs of dimensions $d$ and $d'$, respectively. Suppose a compatible pair $(\psi_1,\psi_2)\colon \PA\times \TA\to \PB\times \TB$ has a lifting
\[
(\psi_1,\tilde{\psi}_2)\colon \PA\times \TMA\to \PB\times \TMB,
\]
where $\tilde{\psi}_2$ is given by
 \[
 \tilde{\psi}_2=\exp\left(
\begin{array}{c c c c}
\ell_{11}  &\ell_{12} &\cdots &\ell_{1m}\\
\ell_{21}  &\ell_{22} &\cdots &\ell_{2m}\\
\vdots  &\vdots &\ddots &\vdots\\
\ell_{m'1}  &\ell_{m'2} &\cdots &\ell_{m'm}
\end{array}\right).
 \]
Let $\{x_1,\ldots,x_{m'}\}$ and $\{y_1,\ldots,y_{m}\}$ be the canonical generators of $\SR(\PB)$ and $\SR(\PA)$. Then for any polynomial $f(x_1,\ldots,x_{m'})\in\SR(\PB)$ we have
\[
\SR(\psi_1,\tilde{\psi}_2)(f(x_1,\ldots,x_{m'}))=f(z_1,\ldots,z_{m'}),
\]
where $z_i=\sum^{m}_{j=1}\ell_{ij}y_j$.    
\end{theorem}

\begin{proof}
Recall from \eqref{dfn_SR[P] isom prototype} that the map $\Phi \colon H_{T^m}(\mathcal{Z}_P) \to \SR(P)$ is induced from the map $\tilde{\varpi}^*$ of \eqref{eq_tilde_varphi}, and similarly for $\Phi'$.
Therefore, the top face of~\eqref{diagram_naturality SR map proof cube} together with~the definition of~$\SR(\psi_1,\tilde\psi_2)$ induces the following commutative diagram
 \[
\begin{tikzcd}
 H^*(B\TMB)\rar{(\tilde{\varpi}')^*}\arrow{d}[left]{(B\tilde{\psi}_2)^*}  &H^*_{\TMB}(\mathcal{Z}_{\PB})\dar{\mathcal{Z}(\psi_1,\tilde{\psi}_2)^*_T}\rar{\Phi'}   &\SR(\PB)\dar{\SR(\psi_1,\tilde{\psi}_2)}\\
 H^*(B\TMA)\rar{\tilde{\varpi}^*}  &H^*_{\TMA}(\mathcal{Z}_{\PA})\rar{\Phi}    &\SR(\PA).
 \end{tikzcd}
\]
Since $(B\tilde{\psi}_2)^*$ is given by 
$(B\tilde{\psi}_2)^*(x_i)=\sum^{m}_{j=1}\ell_{ij}y_j$,  
the statement follows from the commutativity of the above diagram.
\end{proof}

\begin{cor}\label{cor_wSR_n_wSR_bar_morphisms_formula}
With the notation in Proposition~\ref{thm_SR_morphisms_formula}, if $f(x_1,\ldots,x_{m'})$ is a polynomial in $\wSR(\PB,\lambdaB)$ then
\[
\begin{array}{c}
\wSR(\psi_1,\psi_2)(f(x_1,\ldots,x_{m'}))=f(z_1,\ldots,z_{m'});\\[8pt]
\overline{\wSR}(\psi_1,\psi_2)([f(x_1,\ldots,x_{m'})])=[f(z_1,\ldots,z_{m'})],
\end{array}
\]
where $z_i=\sum^{m}_{j=1}\ell_{ij}y_j$. 
Here  $[f(x_1,\ldots,x_{m'})]$ and $[f(z_1,\ldots,z_{m'})]$ are the quotient images of~$f(x_1,\ldots,x_{m'})$ and $f(z_1,\ldots,z_{m'})$ in $\overline{\wSR}(\PB,\lambdaB)$ and $\overline{\wSR}(\PA,\lambdaA)$, respectively.
\end{cor}

\begin{proof}
The corollary follows from Diagrams~\eqref{diagram_naturality of wSR n wSR bar} and~\eqref{diagram_naturality SR map cube} and Proposition~\ref{thm_SR_morphisms_formula}.
\end{proof}

We emphasize that the isomorphism $\overline{\Phi}_X\colon H^*(X(P,\lambda))\to\overline{\wSR}(P,\lambda)$ in~\eqref{eqn_wSR bar isom} depends on the choice of $(P,\lambda)$, while a toric orbifold may be constructed from two different characteristic pairs. For instance, two characteristic pairs $(P, \lambda)$ and~$(P, \lambda')$ satisfying 
\begin{equation}\label{eq_lambda_prime}
\lambda(E)=\epsilon_E\lambda'(E)
\end{equation}
for each facet $E$ of $P$ with  $\epsilon_E\in\{1,-1\}$ define the same toric orbifold by Definition \ref{dfn_toric_orb}. Hence, their corresponding weighted Stanley--Reisner rings $\wSR(P,\lambda)$ and~$\wSR(P, \lambda')$ are isomorphic, though 
they are not the same subrings in $\SR(P)$. The following lemma captures this subtly.

\begin{lemma}\label{lem wSR isom for P with -lambda}
Let $id_P\colon P\to P$ and $id_T\colon T^d\to T^d$ be identity maps. Then~$(id_P,id_T)$ is a compatible pair and
its induced ring morphism
\[
\overline{\wSR}(id_P,id_T)\colon\overline{\wSR}(P,\lambda')\to\overline{\wSR}(P,\lambda)
\]
for characteristic functions $\lambda$ and $\lambda'$ on $P$ satisfying \eqref{eq_lambda_prime}
fits into the commutative diagram
\begin{equation}\label{dgrm_wSR isom for P with -lambda diagram}
\begin{tikzcd}
H^*(X(P,\lambda'))\arrow{rr}{=}\dar{\overline{\Phi}_X'} &&H^*(X(P,\lambda))\dar{\overline{\Phi}_X}\\
\overline{\wSR}(P,\lambda')\arrow{rr}{\overline{\wSR}(id_P,id_T)} &&\overline{\wSR}(P,\lambda)
\end{tikzcd}
\end{equation}
where $\overline{\Phi}_X$ and $\overline{\Phi}'_X$ are the corresponding isomorphisms in~\eqref{eqn_wSR bar isom}. 
Furthermore, if~$f(x_1,\ldots,x_m)$ is a polynomial in $\wSR(P,\lambda')$ then
\[\begin{array}{c}
\wSR(id_P,id_T)(f(x_1,\ldots,x_m))=f(z_1,\ldots,z_m);\\[8pt]
\overline{\wSR}(id_P,id_T)([f(x_1,\ldots,x_{m})])=[f(z_1,\ldots,z_m)],
\end{array}\]
where $z_i=\epsilon_ix_i$ for $1\leq i\leq m$. 
\end{lemma}

\begin{proof}
It is straightforward to check that $(id_P,id_T)$ is a compatible pair in Definition~\ref{def_toric_morphism}. The commutative diagram follows from the bottom face of~\eqref{diagram_naturality of wSR n wSR bar} and the fact that $X(id_P,id_T)\colon X(P,\lambda)\to X(P,\lambda')$ is the identity map.

To show the second part, notice that the compatible pair $(id_P,id_T)$ has a lifting 
\[
(id_P,\tilde{id}_T)\colon P\times T^m\to P\times T^m,
\]
where $\tilde{id}_T\colon T^m\to T^m$ is the Lie group homomorphism
\[
\exp\begin{pmatrix}
\epsilon_1  &   &\\
            &\ddots &\\
            &       &\epsilon_m
\end{pmatrix}.
\]
So Proposition \ref{thm_SR_morphisms_formula} and Corollary~\ref{cor_wSR_n_wSR_bar_morphisms_formula} imply the asserted formulae.
\end{proof}

\begin{eg}\label{eg_wSR isom for triangle with -lambda}
Recall Example~\ref{eg_Xab} and use its notation. Consider $X(\triangle,\lambda_{(a,b)})$ and $X(\triangle,\lambda_{(-a,-b)})$. The ring isomorphism~\eqref{eqn_wSR bar isom} gives two isomorphisms
\[\begin{split}
\overline{\Phi}&\colon H^*(X(\triangle,\lambda_{(a,b)}))\to\overline{\wSR}(\triangle,\lambda_{(a,b)});\\
\overline{\Phi'}&\colon H^*(X(\triangle,\lambda_{(-a,-b)}))\to\overline{\wSR}(\triangle,\lambda_{(-a,-b)}).
\end{split}
\]
Let $[aby_1]\in\overline{\wSR}(\triangle,\lambda_{(a,b)})$ and $[aby_1]'\in\overline{\wSR}(\triangle,\lambda_{(-a,-b)})$ be the equivalence classes of $aby_1$ respectively. They generate the degree 2 components of $\overline{\wSR}(\triangle,\lambda_{(a,b)})$ and~$\overline{\wSR}(\triangle,\lambda_{(-a,-b)})$. Hence their inverse images
\[
\overline{\Phi}^{-1}([aby_1])\quad
\text{and}\quad
(\overline{\Phi'})^{-1}([aby_1]')
\]
generate $H^2(X(\triangle,\lambda_{(a,b)}))=H^2(X(\triangle,\lambda_{(-a,-b)}))$.

On the other hand, Lemma~\ref{lem wSR isom for P with -lambda} implies that
\[
\begin{split}
\overline{\Phi}\left((\overline{\Phi'})^{-1}([aby_1]')\right)
&=\left(\overline{\wSR}(id_{\triangle},id_T)\circ\overline{\Phi}'\right) \left((\overline{\Phi'})^{-1}([aby_1]')\right)\\
&=\overline{\wSR}(id_{\triangle},id_T)([aby_1]')\\
&=-[aby_1].
\end{split}
\]
Hence $\overline{\Phi}^{-1}([aby_1])=-(\overline{\Phi'})^{-1}([aby_1]')$ in $H^2(X(\triangle,\lambda_{(a,b)}))=H^2(X(\triangle,\lambda_{(-a,-b)}))$.
\end{eg}

\section{Toric morphisms of 4-dimensional toric orbifolds}\label{sec_toric morphism in dim four}
From now on, we focus on $2$-dimensional characteristic pairs  $(P, \lambda)$. We consider~$P$ to be an $(n+2)$-gon for some $n\in \mathbb{N}$ and label its vertices and edges by $v_i$ and $E_i$, respectively, for $i=1, \dots, n+2$. For each $E_i$, we write the characteristic vector as~$\lambda(E_i)=(a_i, b_i)\in \mathbb{Z}^2$ as described in Figure~\ref{fig_main char pair}.

\begin{dfn}\label{def_smooth vertex}
A vertex $v_i\in P$ is a \emph{smooth vertex} of $(P,\lambda)$ if
$\{\lambda(E_{i-1}), \lambda(E_{i})\}$ forms an integral basis of $\mathbb{Z}^2$ for $1\leq i\leq n+2$,
where $E_{i-1}$ means $E_{n+2}$ for $i=1$.
\end{dfn}

We note that a smooth vertex $v_i$ of $(P,\lambda)$ corresponds to the fixed point $\pi^{-1}(v_i)$ with the trivial local group in its orbifold chart, hence it is a smooth fixed point in~$X(P, \lambda)$. The readers are referred to \cite[Section 2]{PS} for more details. 
When the context is clear, we simply call $v_{i}$ a smooth vertex of the toric orbifold $X(P,\lambda)$.

When $(P, \lambda)$ has a smooth vertex, we may assume that $v_{n+2}$ is a smooth vertex by relabeling indices and that $\lambda(E_{n+1})=(1,0)$ and $\lambda(E_{n+2})=(0,1)$ by an appropriate basis change of $\mathbb{Z}^2$. In this case, we write $(\underline{a},\underline{b})=\{(a_i,b_i)\}^n_{i=1}$, the first $n$ characteristic vectors of $(P, \lambda)$, and denote by $\lambda_{(\underline{a},\underline{b})}$ the corresponding characteristic function, namely it is defined as  
\begin{equation}\label{eq_1001_char_fun}
\lambda_{(\underline{a},\underline{b})}(E_i)=
\begin{cases}
(a_i, b_i)& i=1, \dots, n;\\
(1,0)& i=n+1;\\
(0,1)& i=n+2,
\end{cases}
\end{equation}
Then, we simply write the associated toric orbifold $X(P, \lambda_{(\underline{a},\underline{b})})$ as 
 $X_{(\underline{a},\underline{b})}$.

The following lemma about the cohomology group of a 4-dimensional toric orbifold $X(P,\lambda)$ is well-known, see for instance \cite[Lemma~3.1]{Fis}, \cite[Theorem 2.5.5]{Jor} and \cite[Corollary~5.1]{KMZ}.

\begin{lemma}\label{lem_evenness_2dim_char_pair}
The cohomology group of a $4$-dimensional toric orbifold $X(P, \lambda)$ is given by 
\begin{equation*}
\begin{tabular}{C{2.4cm}|C{1cm}|C{1cm}|C{1cm}|C{1cm}|C{1cm}|C{1cm}}
$i$		&$0$	&$1$	&$2$	&$3$	&$4$	&$\geq5$\\
\hline
$H^i(X(P,\lambda))$	&$\Z$	&$0$	&$\Z^n$	&$\Z_k$	&$\Z$	&$0$
\end{tabular}
\end{equation*}
where $k=\gcd \left\{|a_ib_j-a_jb_i| \mid 1\leq i<j\leq n+2   \right\}$ and $\Z_k$ means the trivial group if~$k=1$.    
\end{lemma}
Thus, a $2$-dimensional characteristic pair $(P,\lambda)$ is even  if and only if $H^3(X(P,\lambda))$ is trivial (see Definition \ref{dfn_wSR}), which is equivalent to the condition
\[
\gcd \left\{|a_ib_j-a_jb_i| \mid 1\leq i<j\leq n+2   \right\} = 1.
\]
by Lemma~\ref{lem_evenness_2dim_char_pair}. In particular, $(P,\lambda)$ is even if it has a smooth vertex by Definition~\ref{def_smooth vertex}. 

In the following two subsections, we construct two special kinds of toric morphisms: rescaling morphisms and edge-contraction morphisms. Assuming $(P, \lambda)$ is even, we apply 
the theory developed in Section \ref{sec_toric_morphism} to these morphisms, which will be instrumental in proving the main theorem in Section~\ref{sec_pf_main_thm}.

\subsection{Rescaling}\label{subsec_rescaling}
Let $(P,\lambda)$ be an even characteristic pair.
Suppose that there exists an~$i\in\{1,\ldots,n+2\}$ such that $\lambda(E_i)=(a_i,b_i)$ satisfies $a_ib_i\neq0$. Then we define a new characteristic function
\[
\lambda'\colon\{E_1,\ldots,E_{n+2}\}\to\Z^2,
\]
called a \emph{rescaling} of $\lambda$ with respect to $\lambda(E_i)$, by taking
\begin{equation}\label{eqn_rescaled characteristic function}
\lambda'(E_j)=\left(\frac{a_j|b_i|}{g_{ij}},\frac{|a_i|b_j}{g_{ij}}\right)
\end{equation}
for $1\leq j\leq n+2$, where 
$g_{ij}=
\text{gcd}(|a_jb_i|,|a_ib_j|)$.
We set $\gcd(a,0)=\gcd(0,a)=a$ for $a>0$ by convention.
Notice that $\lambda'(E_j)$'s are always primitive vectors. In particular $\lambda'(E_i)=(\epsilon_a, \epsilon_b)$ for the signs $\epsilon_a$ and $\epsilon_b$ of $a_i$ and $b_i$, respectively.

Let $\sigma_i\colon T^2\to T^2$ be defined by
\begin{equation}\label{eq_rescale_torus_map}
\sigma_i(t_1,t_2)=(t_1^{|b_i|},t_2^{|a_i|}).
\end{equation}
Observe that $(id_P,\sigma_i)\colon P\times T^2\to P\times T^2$ is a compatible pair.

\begin{dfn}\label{def_rescaling}
We call the toric morphism
\[
X(id_P,\sigma_i)\colon X(P,\lambda)\to X(P,\lambda')
\]
induced from $(id_P,\sigma_i)$
a \emph{rescaling morphism} of $X(P,\lambda)$ with respect to $\lambda(E_i)$.
\end{dfn}

\begin{lemma}\label{lemma_rescaling morphism wSR}
Let $(P, \lambda)$ be an even characteristic pair such that $\lambda(E_i)=(a_i,b_i)$ satisfies $a_ib_i\neq 0$ for some $i\in \{1,\dots, n\}$, and let $\sigma_i\colon T^2 \to T^2$ be the map defined in \eqref{eq_rescale_torus_map}. For the canonical generators $\{x_1,\ldots,x_{n+2}\}$ of $\SR(P)$, if $f(x_1,\ldots,x_{n+2})$ is a polynomial in $\wSR(P,\lambda')$, then
\[
\overline{\wSR}(id_P,\sigma_i)([f(x_1,\ldots,x_{n+2})])=[f(g_{i1}x_1,\ldots,g_{i,n+2}x_{n+2})],
\]
where $g_{ij}$ is defined in~\eqref{eqn_rescaled characteristic function}.
\end{lemma}

\begin{proof}
The compatible pair $(id_P,\sigma_i)$ has a lifting
\[
(id_P,\tilde{\sigma}_i)\colon P\times T^{n+2}\to P\times T^{n+2}
\]
where $\tilde{\sigma}_i\colon T^{n+2}\to T^{n+2}$ is the Lie group homomorphism given by
\[
\tilde{\sigma}_i(t_1,\ldots,t_{n+2})=(t_1^{g_{i1}},\ldots,t_{n+2}^{g_{i,n+2}}).
\]
One can check that it satisfies Condition~\eqref{dfn_lifting of compatible pair subgp condition} of Definition~\ref{dfn_lifting of compatible pair} by direct computation. Condition~\eqref{dfn_lifting of compatible pair comm diagram} holds since there is a commutative diagram
\[
\begin{tikzcd}
T^{n+2}\rar{\tilde{\sigma}_i}\dar{\exp\Lambda}  &T^{n+2}\dar{\exp\Lambda'}\\
T^2\rar{\sigma_i}                               &T^2
\end{tikzcd}
\]
where $\Lambda$ and $\Lambda'$ are the associated matrices of $\lambda$ and $\lambda'$. Then the lemma follows from Corollary~\ref{cor_wSR_n_wSR_bar_morphisms_formula}.
\end{proof}

\begin{eg}\label{eg_rescaling triangle}
Consider the characteristic pair $(\triangle, \lambda_{(a,b)})$ discussed in Examples~\ref{eg_Xab} and \ref{eg_wSR isom for triangle with -lambda}. 
For the associated toric orbifolds $X_{(a,b)}$ and $X_{(\varepsilon_1,\varepsilon_2)}$ with $\varepsilon_1$ and $\varepsilon_2$ being the signs of $a$ and $b$,
the rescaling morphism
\[
X(id_\triangle,\sigma_1)\colon X_{(a,b)}\to X_{(\varepsilon_1,\varepsilon_2)}
\]
induces the ring homomorphism
\[
\overline{\wSR}(id_\triangle,\sigma_1)\colon \overline{\wSR}(\triangle,\lambda_{(\varepsilon_1,\varepsilon_2)})\to\overline{\wSR}(\triangle,\lambda_{(a,b)}).
\]
Let $\{x_1,x_2,x_3\}$ be the canonical generators of $\SR(\triangle)$. Then we have
\[
\overline{\wSR}(\triangle,\lambda_{(\varepsilon_1,\varepsilon_2)})\cong\Z[x_1,x_2,x_3]/\langle x_1x_2x_3, \varepsilon_1x_1+x_2, \varepsilon_2x_1+x_3\rangle.
\]
Applying Lemma~\ref{lemma_rescaling morphism wSR}, we obtain
\begin{align*}
\overline{\wSR}(id_\triangle,\sigma_1)([\varepsilon_1\varepsilon_2x_1])&= [abx_1];\\
\overline{\wSR}(id_\triangle,\sigma_1)([x_2x_3])&=[|ab|x_2x_3].
\end{align*}
\end{eg}

\subsection{Edge-contraction}\label{section_edge contraction tor orb}
Let $P$ be an $(n+2)$-gon with and $P'$ an $(n'+2)$-gon. Writing $\{v_1,\ldots,v_{n+2}\}$ and $\{v'_1,\ldots,v'_{n'+2}\}$ as vertices of $P$ and $P'$, respectively, their edges are
\begin{align*}
E_i&=\{(1-t)v_i+tv_{i+1}\mid0\leq t\leq 1\} \quad \text{for }1\leq i\leq n+2;\\
E'_j&=\{(1-t)v'_j+tv'_{j+1}\mid0\leq t\leq 1\}   \quad \text{for }1\leq j\leq n'+2.
\end{align*}
We regard $P$ and $P'$ as the cones of their boundaries 
$\partial P$ and~$\partial P'$, respectively, 
that is
$P\cong Cone(\partial P)$
and $P'\cong Cone(\partial P')$.

\begin{dfn}\label{def_edge_contration}
Let $\rho\colon\{1,\ldots,n+2\}\to\{1,\ldots,n'+2\}$ be a surjection that preserves the order. Then the associated \emph{edge-contraction} $\rho\colon P\to P'$ is a continuous map defined as follows.
\begin{itemize}
\item
Its restriction $\partial\rho\colon\partial P\to\partial P'$ is given by
\[
\partial\rho((1-t)v_j+tv_{j+1})=
(1-t)v'_{\rho(j)}+tv'_{\rho(j+1)}
\]
for $0\leq t\leq 1$ and $1\leq k\leq n'+2$, where $v_{n+3}=v_1$ and $v'_{n'+3}=v'_1$;
\item
Define $\rho\colon P\to P'$ to be the cone construction $Cone(\partial\rho)$.
\end{itemize}
\end{dfn}

If $\rho(i)$ and $\rho(i+1)$ are different, then $\rho$ sends the edge $E_i$ homeomorphically onto $E'_{\rho(i)}$. Otherwise $\rho$ contracts $E_i$ to $v'_{\rho(i)}$. For instance, when $P$ is a hexagon and $\rho\colon\{1,\ldots,6\}\to\{1,2,3\}$ is the order-preserving surjection given by
\[
\rho(1)=\rho(2)=1,\quad
\rho(3)=\rho(4)=\rho(5)=2,\quad
\rho(6)=3,
\]
the associated edge-contraction $\rho\colon P\to\triangle$ contracts~$E_1$ to $v_1'$ and $E_3\cup  E_4$ to $v_2'$, which is illustrated in Figure \ref{fig_ex_edge_cont}. 

Given a characteristic pair $(P, \lambda)$ and an edge-contraction $\rho\colon P\to P'$ as in Definition~\ref{def_edge_contration}, one can define a function 
\begin{equation*}
\rho_*\lambda\colon \{E_1', \dots, E'_{n'+2}\} \to \mathbb{Z}^2
\end{equation*}
by $\rho_*\lambda(E'_j)=\lambda(E_i)$ if $\rho(E_i)=E_j'$ for $j=1, \ldots, n'+2$. If $(P', \rho_*\lambda)$ is a characteristic pair, namely it satisfies the conditions of Definition~\ref{dfn_char_funct}, then the map 
\begin{equation}\label{eq_comp_pair_edge_cont}
(\rho, id_T) \colon P\times T^2 \to P'\times T^2
\end{equation}
forms a compatible pair. 
\begin{dfn}\label{dfn_edge contracting morphism 1}
We call the toric morphism  
\[
X(\rho, id_T) \colon X(P, \lambda) \to X(P', \rho_*\lambda)
\]
induced from \eqref{eq_comp_pair_edge_cont} the \emph{edge-contraction morphism} of $X(P,\lambda)$ with respect to $\rho$. 
\end{dfn}

\begin{figure}
\begin{tikzpicture}[scale=0.6]

\draw (30:2)--(90:2)--(150:2)--(210:2)--(270:2)--(330:2)--cycle;

\draw[blue, very thick] (330:2)--(30:2)--(90:2);
\draw[red, very thick] (270:2)--(210:2);

\node at (30:2.5) {\footnotesize$v_4$};
\node at (90:2.5) {\footnotesize$v_5$};
\node at (150:2.5) {\footnotesize$v_6$};
\node at (210:2.5) {\footnotesize$v_1$};
\node at (270:2.5) {\footnotesize$v_2$};
\node at (330:2.5) {\footnotesize$v_3$};

\node at (60:2.2) {\footnotesize$E_4$};
\node at (0:2.2) {\footnotesize$E_3$};
\node at (300:2.2) {\footnotesize$E_2$};
\node at (240:2.2) {\footnotesize$E_1$};
\node at (180:2.2) {\footnotesize$E_6$};
\node at (120:2.2) {\footnotesize$E_5$};

\foreach \a in {30,90,...,330} { 
\draw[fill] (\a:2) circle (2pt); 
}

\draw[fill, blue] (90:2) circle (3pt); 
\draw[fill, blue] (30:2) circle (3pt); 
\draw[fill, blue] (330:2) circle (3pt); 

\draw[fill, red] (270:2) circle (3pt); 
\draw[fill, red] (210:2) circle (3pt);

\node at (0,0) {$\partial P$};

\draw[->] (3.5,0)--(6.5,0);
\node[above] at (5,0) {$\rho$};

\begin{scope}[xshift=350, yshift=0]
\node at (0,0) {$\partial P'$};

\draw (60:2)--(180:2)--(300:2)--cycle;

\foreach \a in {60,180,300} { 
\draw[fill] (\a:2) circle (2pt); 
}
\node at (60:2.5) {\footnotesize$v_2'$};
\node at (180:2.5) {\footnotesize$v_3'$};
\node at (300:2.5) {\footnotesize$v_1'$};

\node at (0:3.2) {\footnotesize$E_1'=\rho(E_2)$};
\node at (150:3) {\footnotesize$E_2'=\rho(E_5)$};
\node at (210:3) {\footnotesize$E_3'=\rho(E_6)$};

\draw[fill,blue] (60:2) circle (3pt); 
\draw[fill,red] (300:2) circle (3pt);

\end{scope}
\end{tikzpicture}
\caption{Example of an edge-contraction.}
\label{fig_ex_edge_cont}
\end{figure}

In the following, we discuss two special types of edge contractions in which $(P,\lambda)$ has a smooth vertex and $P'$ is either a triangle or a square.
We assume that $v_{n+2}$ is a smooth vertex of $(P, \lambda)$ and  $\lambda=\lambda_{(\underline{a},\underline{b})}$ as in~\eqref{eq_1001_char_fun}.

\subsubsection{Edge contractions to triangles}\label{subsubsec_edge_cont_tri}
Suppose $\lambda=\lambda_{(\underline{a},\underline{b})}$ satisfies $a_ib_i\neq 0$ for some~$i\in \{1, \dots, n\}$.
Let $\rho_i\colon P\to\triangle$ be the edge-contraction associated to the order-preserving surjection $\rho_i\colon\{1,\ldots,n+2\}\to\{1,2,3\}$ given by
\begin{equation}\label{eqn_rho_i}
\rho_i(j)=\begin{cases}
1   &\text{if }1\leq j\leq i;\\
2   &\text{if }i+1\leq j\leq n+1;\\
3   &\text{if }j=n+2.
\end{cases}
\end{equation}
Then it contracts all edges of $P$ except $E_i,E_{n+1}$ and $E_{n+2}$. For instance, the edge-contraction of Figure \ref{fig_ex_edge_cont} is $\rho_2$ following the convention of \eqref{eqn_rho_i}. Therefore, we have a compatible pair 
\begin{equation}\label{eq_compatible_pair_rho_i}
    (\rho_i, id_T)\colon P\times T^2 \to \triangle \times T^2
\end{equation} 
with respect to $(P, \lambda)$ and $(\triangle, \lambda_{(a_i, b_i)})$, where the latter is the characteristic pair discussed in Example \ref{eg_Xab}. 
In this case, one can calculate the morphism $\overline{\wSR}(\rho_i, id)$ as in the following Lemma.

\begin{lemma}\label{lemma_edge contract wSR}
Let $(\rho_i, id_T)$ be the compatible pair as in \eqref{eq_compatible_pair_rho_i} with $a_ib_i\neq 0$. For canonical generators $\{x_1, x_2, x_3\}$ and $\{y_1, \dots, y_{n+2}\}$ of $\SR(\triangle)$ and $\SR(P)$,  we have 
 \begin{align*}
\overline{\wSR}(\rho_i,id_T)([a_ib_ix_1])
&=\left[\sum^{i-1}_{k=1}a_kb_iy_k+a_ib_iy_i+\sum^n_{k=i+1}a_ib_ky_k\right],\\
\overline{\wSR}(\rho_i,id_T)([x_2x_3])
&=[y_{n+1}y_{n+2}].
\end{align*}
\end{lemma}

\begin{proof}
First we consider the special case where $a_i=\pm1$ and $b_i=\pm1$. Notice that the compatible pair $(\rho_i,id_T)$ has a lifting
\[
(\rho_i,\tilde{id})\colon P\times T^{n+2}\to\triangle\times T^3
\]
where $\tilde{id}\colon T^{n+2}\to T^3$ is the Lie group homomorphism given by the matrix 
\[
\begin{pmatrix}
a_i a_1     &\cdots &a_ia_{i-1}            &1  &b_ib_{i+1}            &\cdots &b_ib_n        &0  &0\\
0       &\cdots &0                  &0  &a_{i+1}-a_ib_ib_{i+1}    &\cdots &a_n-a_ib_ib_n    &1  &0\\
b_1-a_ib_ia_1 &\cdots &b_{i-1}-a_ib_ia_{i-1}    &0  &0                  &\cdots &0          &0  &1
\end{pmatrix}.
\]
Indeed, one can check Condition~\eqref{dfn_lifting of compatible pair subgp condition} in Definition~\ref{dfn_lifting of compatible pair} by direct computation. Condition~\eqref{dfn_lifting of compatible pair comm diagram} holds since there is a commutative diagram
\[
\begin{tikzcd}
T^{n+2}\rar{\tilde{id}}\dar[swap]{\exp\Lambda}  &T^3\dar{\exp\Lambda'}\\
T^2\rar{id}                               &T^2
\end{tikzcd}
\]
where $\Lambda$ and $\Lambda'$ are characteristic matrices corresponding to $\lambda$ and $\lambda'$, namely
\[
\Lambda=\left(
\begin{array}{c c c c c c c c c}
a_1 &\cdots &a_{i-1}  &a_i  &a_{i+1}    &\cdots &a_n    &1  &0\\
b_1 &\cdots &b_{i-1}  &b_i  &b_{i+1}    &\cdots &b_n    &0  &1
\end{array}
\right)\quad
\text{and}\quad
\Lambda'=\left(
\begin{array}{c c c}
a_i    &1  &0\\
b_i    &0  &1
\end{array}
\right).
\]
Applying Corollary~\ref{cor_wSR_n_wSR_bar_morphisms_formula}, we
obtain
\[
\overline{\wSR}(\rho_i, id_T)([a_ib_ix_1])=\left[\sum^{i-1}_{k=1}a_kb_iy_k+a_ib_iy_i+\sum^n_{k=i+1}a_ib_ky_k\right]
\]
and
\begin{multline}\label{eq_wsr_x2x3}
\overline{\wSR}(\rho_i, id_T)([x_2x_3])\\
=\left[\left(\sum^n_{k=i+1}(a_k-a_ib_ib_k)y_k+y_{n+1}\right)\cdot\left(\sum^{i-1}_{l=1}(b_l-a_ib_ia_l)y_l+y_{n+2}\right)\right].
\end{multline}
Since $y_ky_l=0$ in $\SR(P)$ whenever $k\not\equiv l\pm1\pmod{n+2}$, all terms in the right-hand side of \eqref{eq_wsr_x2x3} vanish except $y_{n+1}y_{n+2}$. Hence we have 
\[
\overline{\wSR}(\rho_i, id_T)([x_2x_3])=[y_{n+1}y_{n+2}]
\]
which proves the claim for $(a_i, b_i)=(\pm1,\pm1)$. 

Next we consider the general case, namely the case where $(a_i, b_i)$ is an arbitrary primitive vector with $a_ib_i\neq 0$. 
The compatible pair $(\rho_i,id_T)\colon P\times T^2\to\triangle\times T^2$ may not have a lifting (see Example \ref{ex_nonexistence_lifting} below), so one cannot apply Corollary~\ref{cor_wSR_n_wSR_bar_morphisms_formula} to compute $\overline{\wSR}(\rho_i,id_T)$ directly. Instead, we consider the commutative diagram
\begin{equation}\label{diagram_edge contract sq to tri n rescaling space}
\begin{tikzcd}
\xq{\underline{a},\underline{b}}\arrow{rr}{{X(\rho_i,id_T)}}\arrow{d}[left]{X(id_{P},\sigma_i)} &&\xabi\arrow{d}{X(id_{\triangle},\sigma_i)}\\
\xq{\underline{a}',\underline{b}’}\arrow{rr}{{X(\rho_i,id_T)'}} &&X_{(\epsilon_1, \epsilon_2)}
\end{tikzcd}
\end{equation}
where 
\begin{itemize}
    \item $id_{P}$ and $id_{\triangle}$ are the identity maps on $P$ and $\triangle$, respectively;
    \item $\sigma_i\colon T^2\to T^2$ is the Lie group homomorphism given by $\sigma_i(t_1,t_2)=(t_1^{|b_i|},t_2^{|a_i|})$;
   \item $\epsilon_1$ and $\epsilon_2$ are the signs of $a_i$ and $b_i$, respectively;
    \item $(\underline{a}', \underline{b}')=\big\{(\frac{a_j|b_i|}{g_{ij}}, \frac{b_j|a_i|}{g_{ij}})\big\}_{j=1}^n$ is the sequence of the characteristic vectors of the rescaling $\lambda'$ of $\lambda$ defined in~\eqref{eqn_rescaled characteristic function};
       \item
    $X(\rho_i,id_T)$ and $X(\rho_i,id_T)'$ are the edge-contractions of $\xq{\underline a,\underline b}$ and $\xq{\underline a',\underline b'}$, respectively, induced by $\rho_i\colon P\to\triangle$ and $id_T\colon T^2\to T^2$.
\end{itemize}  
It induces the following cubical diagram whose faces commute
\begin{equation}\label{eq_cubical_diag}
\small
\begin{tikzcd}
&H^*(X_{(\epsilon_1, \epsilon_2)})
\ar{rr}{(X(\rho_i,id_T)')^*}
\ar{dd}{}
\ar{dl}[swap]{X(id_{\triangle},\sigma_i)^*}
  &   &H^*(\xq{\underline a', \underline b'})
\ar{dd}\ar{dl}{X(id_P,\sigma_i)^*}\\
H^*(\xabi)
\ar[crossing over,near end]{rr}{X(\rho_i,id_T)^*}\ar{dd}   &   &
H^*(\xq{\underline a, \underline b})   &\\
    &\overline{\wSR}(\triangle,\lambda_{(\epsilon_1,\epsilon_2)})
\ar{dl}[swap,near start]{\overline{\wSR}(id_{\triangle},\sigma_i)}
\ar{rr}[near start]{\overline{\wSR}(\rho_i,id_T)'}   &   &
\overline{\wSR}(P,\lambda')
\ar{dl}{\overline{\wSR}(id_P,\sigma_i)}\\
\overline{\wSR}(\triangle,\lambda_{(a_i,b_i)})
\ar{rr}{\overline{\wSR}(\rho_i,id_T)}   &    &\overline{\wSR}(P,\lambda),\ar[crossing over, leftarrow]{uu}
&
\end{tikzcd}
\end{equation}
where all vertical maps are the isomorphism studied in~\eqref{eqn_wSR bar isom}. The top face is induced by~\eqref{diagram_edge contract sq to tri n rescaling space}. The front, the rear, the left and right faces commute due to the definitions of $\overline{\wSR}(\rho_i,id_T),\overline{\wSR}(\rho_i,id_T)',\overline{\wSR}(id_{\triangle},\sigma_i)$ and $\overline{\wSR}(id_P,\sigma_i)$, respectively. Then a diagram chasing shows that the bottom face commutes as well.

Recall from Lemma \ref{lemma_rescaling morphism wSR} that 
\[
\overline{\wSR}(id_\triangle, \sigma_i)([\epsilon_1\epsilon_2x_1])=[a_ib_ix_1].
\]
Hence, the commutativity of the bottom face of \eqref{eq_cubical_diag} implies that 
\begin{align*}
\overline{\wSR}(\rho_i,id_T)([a_ib_ix_1]) &=\left(\overline{\wSR}(\rho_i,id_T)\circ\overline{\wSR}(id_{\triangle},\sigma_i)\right)([\epsilon_1 \epsilon_2 x_1])\\
&=\left( \overline{\wSR}(id_P,\sigma_i)\circ\overline{\wSR}(\rho_i,id_T)'\right)([ \epsilon_1 \epsilon_2 x_1])\\
&=\overline{\wSR}(id_P,\sigma_i)\left(\left[\sum^{i-1}_{k=1} \epsilon_2 a'_ky_k+ \epsilon_1 \epsilon_2  y_i+\sum^n_{k=i+1}\epsilon_1 b'_ky_k\right]\right)\\
&=\left[\sum^{i-1}_{k=1}\epsilon_2 a'_k g_{ik} y_k+ \epsilon_1 \epsilon_2 g_{ii}  y_i+\sum^n_{k=i+1}  \epsilon_1b'_k g_{ik} y_k\right]\\
&=\left[\sum^{i-1}_{k=1}a_kb_iy_k+a_ib_iy_i+\sum^n_{k=i+1}a_ib_ky_k\right],
\end{align*}
where the last equality follows because $g_{ik}=\gcd(|a_kb_i|, |b_ka_i|)$,  $a'_k=\frac{a_k|b_i|}{g_{ik}}$ and $b'_k=\frac{b_k|a_i|}{g_{ik}}$.
A similar calculation shows that $\overline{\wSR}(\rho_i,id_T)([x_2x_3])=[y_{n+1}y_{n+2}]$.
\end{proof}

Here we provide an example of an
edge-contraction that does not have a lifting. It shows the reason why Corollary~\ref{cor_wSR_n_wSR_bar_morphisms_formula} cannot be applied directly to the general case in the proof of Lemma~\ref{lemma_edge contract wSR}.
 
\begin{eg}\label{ex_nonexistence_lifting}

Let $(\square,\lambda)$ and $(\triangle,\lambda')$ be characteristic pairs where $\lambda$ and $\lambda'$ have associated characteristic matrices
\[\Lambda=\begin{pmatrix} 2 & -3 & 1 & 0 \\1 & -2 & 0 & 1\end{pmatrix}\quad
\text{and} \quad
\Lambda'=\begin{pmatrix} -3 & 1 & 0 \\ -2 & 0 & 1\end{pmatrix}
\]
respectively.
Consider $\rho_2\colon \square\to \triangle$ in~\eqref{eqn_rho_i} and the induced edge-contraction morphism 
\[
X(\rho_2, id_T) \colon \xq{2,1),(-3,-2}\to \xq{-3,-2}.
\]
We claim that $X(\rho_2, id_T)$ does not admit a lifting. 
Assume that there exists a lifting $\tilde{id}_T\colon T^4 \to T^3$ making the following diagram commute
\begin{equation}\label{eq_lifting_id}
\begin{tikzcd}
T^4\arrow{r}{\tilde{id}_T} \arrow{d}[swap]{\exp \Lambda} &T^3 \arrow{d}{\exp\Lambda'}\\
T^2\arrow{r}{=} &T^2.
\end{tikzcd}
\end{equation}
The first condition of Definition \ref{dfn_lifting of compatible pair} implies that $\tilde{id}_T(t,1,1,1)=(t^{a_1}, {1, t^{a_2}})$ for some~$a_1, a_2 \in \mathbb{Z}$ as $\rho_2(E_1)=E_1'\cap E_3'$. Observe that 
\[
(\exp \Lambda'\circ \tilde{id}_T) (t,1,1,1)=(\exp \Lambda') (t^{a_1},  1,t^{a_2})=(t^{-3a_1}, t^{-2a_1+a_2}), 
\]
which has to agree with 
\[
(\exp \Lambda)(t,1,1,1)=(t^2,t)
\]
by the commutativity of \eqref{eq_lifting_id}. This contradicts to the fact that $a_1$ and $a_2$  are integers. 
\end{eg}

\subsubsection{Edge contractions to squares}\label{subsubsec_edge_cont_sq}
Suppose $\lambda=\lambda_{(\underline{a},\underline{b})}$ in~\eqref{eq_1001_char_fun} satisfying
\[
a_i\neq0,\quad
b_j\neq0,
\quad\text{and}\quad
a_ib_j-a_jb_i\neq0.
\]
for some $i,j\in\{1,\ldots,n\}$ with $i<j$.
Let $\square$ denote a square and let $\lambda_{ij}$ be a characteristic function on $\square$ with the characteristic matrix 
\[
\left(\begin{array}{c c c c}
a_i &a_j    &1  &0\\
b_i &b_j    &0  &1
\end{array}\right).
\]
Let $\rho_{ij}\colon P\to\square$ be the edge-contraction associated with the order-preserving surjection $\rho_{ij}\colon\{1,\ldots,n+2\}\to\{1,2,3,4\}$ given by
\begin{equation}\label{eq_rho_ij}
\rho_{ij}(k)=\begin{cases}
1   &\text{if }1\leq k\leq i;\\
2   &\text{if }i+1\leq k\leq j;\\
3   &\text{if }j+1\leq k\leq n+1;\\
4   &\text{if }k=n+2.
\end{cases}
\end{equation}
Then, the edge-contraction $\rho_{ij}$ together with the identity map $id_T$ on $T^2$ gives us 
a compatible pair $(\rho_{ij},id_T)\colon P\times T^2\to \square \times T^2$.

\begin{lemma}\label{lemma_rho_ij to sq}
 Let $\{x_1,\ldots,x_4\}$ and $\{y_1,\ldots,y_{n+2}\}$ be the canonical generators of~$\SR(\square)$ and $\SR(P)$, respectively.
If $\lambda(E_i)=(\pm1,0)$ and $\lambda(E_j)=(a_j,\pm1)$ for some $a_j\in \mathbb{Z}$ then
\[
\overline{\wSR}(\rho_{ij},id)([f(x_1,\cdots,x_4)])=f(z_1,\ldots,z_4)
\]
where
\begin{align*}
z_1&=a_1a_iy_1+\cdots+a_{i-1}a_iy_{i-1}+y_i+a_i(a_{i+1}-a_jb_jb_{i+1})y_{i+1}+\cdots+\\
&\quad\,\, a_i(a_{j-1}-a_jb_jb_{j-1})y_{j-1}\\
z_2&=b_ja_{i+1}y_{i+1}+\cdots+b_ja_{j-1}y_{j-1}+y_j+b_jb_{j+1}y_{j+1}+\cdots+b_jb_ny_n\\
z_3&=(a_{j+1}-a_jb_jb_{j+1})y_{j+1}+\ldots+(a_{n}-a_jb_jb_{n})y_{n}+y_{n+1}\\
z_4&=b_1y_1+\ldots+b_{i-1}y_{i-1}+y_{n+2}.
\end{align*}
\end{lemma}

\begin{proof}
We verify the claim by calculating the lifting  $(\rho_{ij},\tilde{id})$ of $(\rho_{ij},id)$ directly. Indeed, one can take $\tilde{id}\colon T^{n+2}\to T^4$ to be the exponential of the map given by the following matrix 
\[
\begin{blockarray}{ccccccc}
\begin{block}{(c|c|c|c|c|cc)}
A_1 & \begin{matrix}
1\\
0\\
0\\
0
\end{matrix} & A_2
&
\begin{matrix}
0\\
1\\
0\\
0
\end{matrix} & A_3 &
\begin{matrix}
0\\
0\\
1\\
0
\end{matrix} &
\begin{matrix}
0~~\\
0~~\\
0~~\\
1~~
\end{matrix}\\
\end{block}
&\small{\overset{\uparrow}{\text{$i$-th}}}&&\small{\overset{\uparrow}{\text{$j$-th}}}&&&
\end{blockarray}
\]
where
\begin{itemize}
\item
$A_1=\begin{pmatrix}
a_1a_i  &\cdots &a_{i-1}a_i\\
0  &\cdots &0\\
0  &\cdots &0\\
b_1  &\cdots &b_{i-1}
\end{pmatrix}$
\smallskip 
\item
$A_2=\begin{pmatrix}
a_ia_{i+1}-a_ia_jb_jb_{i+1}  &\cdots &a_ia_{j-1}-a_ia_jb_jb_{j-1}\\
b_ja_{i+1}  &\cdots &b_ja_{j-1}\\
0  &\cdots &0\\
0  &\cdots &0
\end{pmatrix}$
\smallskip 
\item
$A_3=\begin{pmatrix}
0  &\cdots &0\\
b_jb_{j+1}  &\cdots &b_jb_{n}\\
a_{j+1}-a_jb_jb_{j+1}  &\cdots &a_{n}-a_jb_jb_{n}\\
0  &\cdots &0
\end{pmatrix}$.
\end{itemize}
Then the lemma follows from Corollary~\ref{cor_wSR_n_wSR_bar_morphisms_formula}.
\end{proof}

\begin{lemma}\label{lem_5_deg_4_are_all_same}
Let $(P,\lambda)$ be a characteristic pair with $\lambda=\lambda_{(\underline{a},\underline{b})}$ in~\eqref{eq_1001_char_fun}, and let~$\{y_1,\ldots,y_{n+2}\}$ be the canonical generators of $\SR(P)$. Then $[y_{n+1}y_{n+2}]$ is a degree-4 generator of $\overline{\wSR}(P,\lambda)$.
\end{lemma}

\begin{proof}
Suppose that there exists an $i\in\{1,\ldots,n\}$ such that $\lambda(E_i)=(a_i,b_i)$ satisfies
\[
a_ib_i\neq0.
\]
Then $(\triangle,\lambda_{(a_i,b_i)})$ is a characteristic pair.
Let $\rho_i\colon P\to \triangle$ be the edge contraction given in~\eqref{eqn_rho_i},
and let $\{z_1,z_2,z_3\}$
be the canonical generator of $\SR(\triangle)$.
From Example~\ref{eg_Xab}, we know $[z_2z_3]$ is a degree-4 generator of $\overline{\wSR}(\triangle,\lambda_{(a_i,b_i)})$.
By Lemma~\ref{lemma_edge contract wSR}, it follows that
$[y_{n+1}y_{n+2}]=\overline{\wSR}(\rho_i,id_T)[z_2z_3]$ is a degree-4 generator of $\overline{\wSR}(P,\lambda)$.

Suppose for $1\leq i\leq n$ all characteristic vectors $\lambda(E_i)=(\pm1,0)$ or $(0,\pm1)$. 
Since~$\lambda$ is a characteristic function,
there must exist $j<k\in\{1,\ldots,n\}$ such that
\[
\lambda(E_j)=(\pm1,0)
\quad\text{and}\quad
\lambda(E_k)=(0,\pm1).
\]
Then $(\square,\lambda_{(a_j,b_j),(a_k,b_k)})$ is a characteristic pair and we have
\[
\overline{\wSR}(\square,\lambda_{(a_j,b_j),(a_k,b_k)})=\SR(\square)/\langle{\pm x_1+x_3,\pm x_2+x_4}\rangle\cong\Z[x_3,x_4]/\langle{x_3^2,x_4^2}\rangle.
\]
Therefore, $[x_3x_4]$ is a degree-4 generator of $\overline{\wSR}(\square,\lambda_{(a_j,b_j),(a_k,b_k)})$.
By Lemma~\ref{lemma_rho_ij to sq}
\[
[y_{n+1}y_{n+2}]=\overline{\wSR}(\rho_{jk},id)^*([x_3x_4])
\]
is a degree-4 generator of $\overline{\wSR}(P,\lambda)$. Hence, the proof is complete.
\end{proof}

\section{Cellular bases}\label{sec_DTS}

In Section~\ref{sec_preliminaries}, we have identified $H^*(X(P,\lambda))$ with the quotient ring $\overline{\wSR}(P,\lambda)$ via Isomorphism~\eqref{eqn_wSR bar isom}.
In this section, we provide an alternative description of $H^*(X(P,\lambda))$ using \emph{cellular bases}, following the approach of~\cite{FSS}. Namely, we define a specific cellular structure for $X(P,\lambda)$ (Example~\ref{eg_CW str of toric orb}) and use it to construct an additive basis for its cohomology (Definition~\ref{def_cellular_basis_cel_cup_repn}). This topological perspective has the advantage of directly associating each cell in $X(P,\lambda)$ with a corresponding cohomology class in $H^*(X(P,\lambda))$.

To connect these two descriptions of $H^*(X(P,\lambda))$, we study the cohomology morphisms induced by the rescaling and edge-contraction morphisms from the previous section. 
For this purpose, we extend $4$-dimensional toric orbifolds to a broader class of spaces, called \emph{degenerate toric spaces}, and generalize toric morphisms accordingly. This section begins with a review of the category of mapping cones and their associated cellular bases. We then define degenerate toric spaces and compute the cohomology morphisms induced by rescaling and edge-contraction morphisms.

\subsection{The category of mapping cones}\label{subsec_cat_mapping_cone}
Let $\mathscr{C}_n$
be the full subcategory of topological spaces that consists of CW-complexes of the form
\begin{equation}\label{eqn_C_f cell str}
C_f=\left(\bigvee^n_{i=1}S^2_i\right)\cup_f D^4
\end{equation}
where $S^2_i$ is a $2$-sphere with label $i$ and $f\colon S^3\to \bigvee_{i=1}^nS^2_i$ is the attaching map of the $4$-cell $D^4$.
Obviously $C_f$ is the mapping cone of $f$ and its homotopy type is determined by the homotopy class of $f$.

\begin{eg}\label{eg_CW str of toric orb}
Following the idea in~\cite[Section 4]{BSS} and~\cite[Section 2]{FSS}, we show that every $4$-dimensional toric orbifold $X(P, \lambda)$ with a smooth vertex 
is homotopy equivalent to a mapping cone $C_f$ in $\mathscr{C}_n$.
Without loss of generality, we assume that~$v_{n+2}$ is a smooth vertex.
Draw a line segment $L$ that intersects $E_{n+1}$ and~$E_{n+2}$ as in Figure~\ref{fig_q CW}.
The neighborhood of~$v_{n+2}$ bounded by $L$ is a small triangle $Cone(L)$ and we have
\begin{equation}\label{eq_Y decomposition}
X(P,\lambda)\cong\pi^{-1}(Cone(L))\cup\pi^{-1}(P\setminus\{v_{n+2}\}),
\end{equation}
where $\pi \colon X(P, \lambda) \to P$ is the orbit map.
On the right hand side, $\pi^{-1}(Cone(L))$ is a $4$-disk since it is the cone of $\pi^{-1}(L)\cong S^3$, and $\pi^{-1}(P\setminus\{v_{n+2}\})$ is homotopy equivalent to $\bigvee^n_{i=1}S^2_i$ via the composite
\[
\pi^{-1}(P\setminus\{v_{n+2}\})\overset{\simeq}{\longrightarrow}\pi^{-1}\left(\bigcup^n_{i=1}E_i\right)\cong\bigcup^n_{i=1}\pi^{-1}(E_i)\overset{\simeq}{\longrightarrow}\bigvee^n_{i=1}S^2_i
\]
where the first map is induced by the deformation retract $P\setminus\{v_{n+2}\}\to\bigcup^n_{i=1}E_i$ and the last map identifies $\pi^{-1}(E_i)$ with a pointed $2$-sphere $S^2_i$ by shrinking the line $E_i\times\{(1,1)\}\subset\pi^{-1}(E_i)$ to a point.
Equation~\eqref{eq_Y decomposition} implies that $X(P,\lambda)$ is the mapping cone of the composite
\begin{equation}\label{eq_f attaching map}
f\colon S^3\cong\pi^{-1}(L)\hookrightarrow\pi^{-1}(P\setminus\{v_{n+2}\})\simeq\bigvee^n_{i=1}S^2_i,
\end{equation}
and hence is a CW-complex in $\mathscr{C}_n$.
\begin{figure}
\begin{tikzpicture}[scale=0.6]

\draw[fill=yellow!50, yellow!50, thick](30:2)--(90:2)--(150:2)--(210:2)--(270:2)--(330:2)--cycle;
\draw[fill=teal!30, teal!30, thick] (60:1.7)--(90:2)--(120:1.7)--cycle;
\draw[red, very thick] (60:1.7)--(120:1.7);
\draw[thick](30:2)--(90:2)--(150:2)--(210:2)--(270:2)--(330:2)--cycle;

\foreach \vangle in {30,90,150,210,270,330}
\draw[fill] (\vangle:2) circle (2pt);

\node[right] at (30:2) {\footnotesize$v_{n+1}$};
\node[above] at (90:2) {\footnotesize$v_{n+2}$};
\node[left] at (150:2) {\footnotesize$v_1$};
\node[left] at (210:2) {\footnotesize$v_2$};
\node[below] at (270:2) {\footnotesize$\cdots$};
\node[right] at (330:2) {\footnotesize$v_{n}$};

\node[right] at (60:2) {\footnotesize$E_{n+1}$};
\node[left] at (120:2) {\footnotesize$E_{n+2}$};
\node[left] at (180:2) {\footnotesize$E_1$};
\node[right] at (0:2) {\footnotesize$E_{n}$};
\node at (90:1) {$L$};
\node at (0:4.2) {$=$};

\begin{scope}[xshift=230]
\draw[fill=yellow!50, thick](30:2)--(60:1.7)--(120:1.7)--(150:2)--(210:2)--(270:2)--(330:2)--cycle;
\draw[red, very thick] (60:1.7)--(120:1.7);
\node at (90:1) {$L$};
\foreach \vangle in {30,150,210,270,330}
\draw[fill] (\vangle:2) circle (2pt); 

\node[right] at (30:2) {\footnotesize$v_{n+1}$};
\node[left] at (150:2) {\footnotesize$v_1$};
\node[left] at (210:2) {\footnotesize$v_2$};
\node[below] at (270:2) {\footnotesize$\cdots$};
\node[right] at (330:2) {\footnotesize$v_{n}$};

\node[left] at (180:2) {\footnotesize$E_1$};
\node[right] at (0:2) {\footnotesize$E_{n}$};

\draw[fill=white] (60:1.7) circle (2pt); 
\draw[fill=white] (120:1.7) circle (2pt);

\node at (0:4.2) {\Large$\cup$};

\end{scope}

\begin{scope}[xshift=410]

\draw[fill=teal!40, thick] (60:1.7)--(90:2)--(120:1.7)--cycle;

\draw[fill] (90:2) circle (2pt);
\node[above] at (90:2) {\footnotesize$v_{n+2}$};


\draw[red, very thick] (60:1.7)--(120:1.7);
\node at (90:1) {$L$};

\draw[fill=white] (60:1.7) circle (2pt); 
\draw[fill=white] (120:1.7) circle (2pt); 

\end{scope}
\end{tikzpicture}
\caption{$X(P, \lambda)$ as a mapping cone.}\label{fig_q CW}
\end{figure}
\end{eg}

\begin{eg}\label{ex_mapping_cone_not_toric_orb}
Let $X(P,\lambda)$ be a $4$-dimensional toric orbifold with a smooth vertex~$v_{n+2}$. Choose an interior point $v\in E_i$ for some $i\in\{1,\ldots,n+2\}$ and let
\[
Y=X(P,\lambda)/\pi^{-1}(v).
\]
Since $\pi^{-1}(v)$ is homeomorphic to $S^1$ and $\pi^{-1}(v)\hookrightarrow X(P,\lambda)$ is null homotopic, we have $
Y\simeq X(P,\lambda)\vee S^2$.
By Example~\ref{eg_CW str of toric orb} $X(P,\lambda)$ is in $\mathscr{C}_n$, so $Y$ is in $\mathscr{C}_{n+1}$.
Note that it is not a toric orbifold.
\end{eg}

For any CW-complex $C_f\in\mathscr{C}_n$, its cohomology is
\[
H^i(C_f)\cong\begin{cases}
\Z^n		&i=2\\
\Z		&i=0,4\\
0		&\text{otherwise},
\end{cases}
\]
whose generators correspond to its $2$-cells $S^2_i$ and $4$-cell $D^4$ in~\eqref{eqn_C_f cell str}.

\begin{dfn}\label{def_cellular_basis_cel_cup_repn}
For an object $C_f\in\mathscr{C}_n$, its \emph{cellular basis}
 \[
 \mathcal{B}=\{u_1,\ldots,u_n;v\}
 \]
 is a basis for $\tilde{H}^*(C_f)$, where $u_i\in H^2(C_f)$ and $v\in H^4(C_f)$ are cohomology classes corresponding to the $2$-sphere $S^2_i$ and the $4$-cell $D^4$.
 In addition, we call the $(n\times n)$-integral matrix
 \[
 \mcup{\mathcal B, C_f} \colonequals  (c_{ij})_{1\leq i,j\leq n},
 \]
 where $c_{ij}$ satisfies $u_i \cup u_j = c_{ij} v$, the \emph{cellular cup product representation} of $C_f$ with respect to $\mathcal{B}$.
\end{dfn}

\begin{remark}
Given a CW-complex $C_f\in\mathscr{C}_n$, its cellular basis is not unique. For example, if $\{u_1,\ldots,u_n;v\}$ is a cellular basis for $C_f$ then $-u_i$ and $-v$ correspond to the same cells $S^2_i$ and $D^4$ but with opposite orientations. Hence
\[
\{ \epsilon_1 u_1, \dots, \epsilon_n u_n ; \epsilon_{n+1} v\}
\]
with $\epsilon_i=\pm 1$ for $i=1, \dots, n+1$
is also a cellular basis of~$C_f$.
\end{remark}

\subsection{Degenerate toric spaces}\label{subsec_dts}

The quotient space $Y$ in Example~\ref{ex_mapping_cone_not_toric_orb} can be constructed in a way similar to the construction of toric orbifolds, except that the characteristic vectors of adjacent edges of $P$ need not be linearly independent.

\begin{dfn}\label{def_deg_char_fun}
Let $P$ be an $(n+2)$-gon as in Figure~\ref{fig_main char pair}. 
\begin{enumerate}
\item
A \emph{degenerate characteristic function} on $P$ is a map
\[
\eta\colon \{E_1,\ldots,E_{n+2}\}\to\Z^2
\]
sending edges of $P$ to primitive vectors $\eta(E_i)=(a_i,b_i)$, and the pair $(P,\eta)$ is called a \emph{degenerate characteristic pair}.
\item The \emph{degenerate toric space} associated with a degenerate characteristic pair $(P,\eta)$ is the quotient space
\[
X(P, \eta)=P\times T^2/_{\sim_{d}}
\]
equipped with a $T^2$-action and a map $\pi\colon X(P, \eta) \to P$ defined as follows:
\begin{itemize}
\item[(i)]
The equivalence relation $(x,t)\sim_d (y,s)$ is given by $x=y$ and
\begin{enumarabic}
\item
$t^{-1}s\in T_{E(x)}$ when $x$ is not a vertex of $P$, or
\item
$t$ and $s$ are any points in $T^2$ when $x$ is a vertex of $P$.
\end{enumarabic}
We denote by $[x,t]_{\sim_d}$ the equivalence class of $(x,t)\in P\times T^2$. 
\item[(ii)]
The $T^2$-action $T^2\times X(P,\eta) \to X(P,\eta)$ is given by
$(g,[x,t]_{\sim_d})\mapsto [x,g\cdot t]_{\sim_d}$ for any $g\in T^2$.
\item[(iii)]
The orbit map $\pi\colon X(P, \eta) \to P$ is given by
$\pi([x,t]_{\sim_d})=x$.
\end{itemize}
\end{enumerate}
\end{dfn}

Comparing Definitions~\ref{dfn_char_funct} and~\ref{def_deg_char_fun}, we see that a degenerate characteristic function $\eta$ becomes a characteristic function if the image vectors of any two adjacent edges are linearly independent. In this case, Condition~(i) of Definition \ref{def_deg_char_fun} (2) agrees with the first condition of Definition \ref{dfn_toric_orb}, and hence $X(P,\eta)$ is a toric orbifold.

Observe that if $\eta$ is a characteristic function of Definition~\ref{dfn_char_funct}, then Condition~(i) of Definition \ref{def_deg_char_fun} (2) agrees with the first condition of Definition \ref{dfn_toric_orb}.

\begin{eg}
In Example~\ref{ex_mapping_cone_not_toric_orb} the space $Y=X(P,\lambda)/\pi^{-1}(v)$ can be regarded as a degenerate toric space $X(P',\eta)$, where $P'$ is an $(n+3)$-gon and $\eta$ is given by
\[
\eta(E_j)=\begin{cases}
\lambda(E_j)		&\text{for }1\leq j\leq i\\
\lambda(E_i)		&\text{for }j=i+1\\
\lambda(E_{j-1})	&\text{for }i+2\leq j\leq n+3.
\end{cases}
\]
Note that $\eta(E_i)=\eta(E_{i+1})$ so $\eta$ is not a characteristic function in Definition~\ref{dfn_char_funct}.
\end{eg}

As in Definition~\ref{def_smooth vertex} we still call a vertex $v_i\in P$ a \emph{smooth vertex} of $(P,\eta)$ if~$\{\eta(E_{i-1}),\eta(E_i)\}$ forms a basis for $\Z^2$ or equivalently it satisfies
\[
|a_{i-1}b_i-a_ib_{i-1}|=1.
\]
When a degenerate toric space $X(P,\eta)$ has a smooth vertex, the argument in Example~\ref{eg_CW str of toric orb} applies. Therefore, $X(P,\eta)$ is the mapping cone of an attaching map $f$ of the form~\eqref{eq_f attaching map} and has a cellular basis.

Moreover, after relabeling vertices of $P$ and applying basis change in $\Z^2$, every degenerate toric space with a smooth vertex is homeomorphic to
$X_{(\underline{a},\underline{b})}=X(P,\eta_{(\underline{a},\underline{b})})$
where $(\underline{a},\underline{b})=\{(a_i,b_i)\}^n_{i=1}$ is a sequence  of primitive integer vectors and
\begin{equation}\label{eq_1001_deg_char_fun}
\eta_{(\underline{a},\underline{b})}(E_i)=\begin{cases}
(a_i,b_i)	&1\leq i\leq n\\
(1,0)		&i=n+1\\
(0,1)		&i=n+2.
\end{cases}
\end{equation}

\begin{eg}\label{eg_X_(1,0)}
We verify that $X_{(1,0)}$ is homotopy equivalent to $S^4 \vee S^2$. Since $P$ is a triangle and $X_{(1,0)}$ has a smooth vertex $v_3$, it is the mapping cone of an attaching map $f\colon S^3 \to S^2$ in the form of~\eqref{eq_f attaching map}. It suffices to show that $f$ is null homotopic.
By~\eqref{eq_f attaching map} it factors through the composite
\[
f'\colon\pi^{-1}(L)\hookrightarrow\pi^{-1}(P\setminus \{v_3\})\xrightarrow{\text{retract}}\pi^{-1}(E_1).
\]
Parameterizing $L$ and $E_1$ by affine maps $L,E_1\colon[0,1]\to P$ such that
\[
L(0)\in E_3,\quad L(1)\in E_2,\quad E_1(0)=v_1\quad \text{and}\quad E_1(1)=v_2,
\]
we may regard $\pi^{-1}(L)$ and $\pi^{-1}(E_1)$ as quotient spaces of~$[0,1]\times T^2$, and write their elements as $[t, x_1, x_2]_L$ and $[t, x_1, x_2]_E$ respectively. Then we have 
\[
f'([t, x_1, x_2]_L)=[t, x_1,x_2]_E.
\]
Since $\eta_{(1,0)}(E_1)=\eta_{(1,0)}(E_2)$, there is a null-homotopy $H\colon \pi^{-1}(L)\times[0,1]\to\pi^{-1}(E_1)$
of $f'$
given by
\[
H([t, x_1, x_2]_L,\theta)=[(1-\theta)t, x_1, x_2]_E.
\]
Therefore $f$ is null-homotopic and $X_{(1,0)}\simeq S^4\vee S^2$.
\end{eg}

\subsection{Toric morphisms of degenerate toric spaces}\label{section_edge contraction revisited}

The definitions of compatible pairs $(\psi_1,\psi_2)\colon P\times T^2\to P\times T^2$ and toric morphisms in Definition~\ref{def_toric_morphism} apply directly to degenerate toric spaces.

In what follows, we assume that degenerate toric spaces have smooth vertices and are of the form $X_{(\underline{a},\underline{b})}$ for some $(\underline{a},\underline{b})$ in~\eqref{eq_1001_deg_char_fun}. Using their cellular bases, we
determine the cohomology morphisms induced by rescaling and edge-contraction morphisms in several special cases needed for later calculations.

\subsubsection{Rescaling}
Consider $X_{(1,0)}$ in Example~\ref{eg_X_(1,0)}. 
For any positive integers $a$ and~$b$, the rescaling morphism (see~\eqref{eq_rescale_torus_map})
\[
X(id,\sigma)\colon X_{(1,0)}\to X_{(1,0)}
\]
induced by the compatible pair $(id,\sigma)\colon \triangle\times T^2\to \triangle\times T^2$ with
$\sigma=\exp\begin{pmatrix}
b	&0\\
0	&a
\end{pmatrix}$
is always well-defined.

\begin{lemma}\label{lemma_rescale Y_(1,0)}
Let $\{1,u_1,v\}$ be a cellular basis for $H^*(X_{(1,0)})$. Then
\[
X(id,\sigma)^*(u_1)=au_1
\quad\text{and}\quad
X(id,\sigma)^*(v)=abv.
\]
\end{lemma}

\begin{proof}
Recall that $u_1$ is the cohomology class corresponding to the 2-sphere $\pi^{-1}(E_1)$, which is the suspension of~$T^2/T_{E_1}\cong\{1\}\times S^1$. Since the restriction of $X(id_{\triangle},\sigma)$ to $\pi^{-1}(E_1)$ has degree $a$, we have $X(id_{\triangle},\sigma)^*(u_1)=au_1$. Similarly, one can show that $X(id_{\triangle},\sigma)^*(v)=abv$.
\end{proof}

\subsubsection{Edge-contraction}
Given a degenerate characteristic pair $(P, \eta)$ and
an edge-contraction $\rho\colon P\to P'$ given by
an order-preserving surjection
\[
\rho\colon\{1,\ldots,n+2\}\to\{1,\ldots,n'+2\},
\]
one can define the edge contraction $X(\rho,id)\colon X(P,\eta)\to Y(P',\rho_*\eta)$ as in Definition~\ref{def_edge_contration}.
Assume that $X(P,\eta)=X_{(\underline{a},\underline{b})}$ and $\rho$ satisfy
\begin{equation}\label{eqn_deg edge contraction condition}
\rho(n+1)=n'+1
\quad\text{and}\quad
\rho(n+2)=n'+2.
\end{equation}
Then $X(P',\rho_*\eta)=X_{(\underline{a}',\underline{b}')}$ where
$(\underline a',\underline b')=\{(a_{\sur_k},b_{\sur_k})\}^{n'}_{k=1}$ with
$$\sur_k=\max\{j\mid\rho(v_j)=v_k\}.$$

\begin{lemma}\label{lemma_rho_S n cellular basis}
Let $\rho\colon P\to P'$ satisfy~\eqref{eqn_deg edge contraction condition}. If
$\{u_1,\ldots,u_n;v\}$ and $\{u'_1,\ldots,u'_{n'};v'\}$ are cellular bases of 
$H^*(X_{(\underline a,\underline b)})$ and $H^*(X_{(\underline a',\underline b')})$ respectively, then
\[
X(\rho,id_T)^*\colon H^*(X_{(\underline a',\underline b')})\to H^*(X_{(\underline a,\underline b)})
\]
sends $v'$ to $\pm v$, and sends $u'_k$ to $\pm u_{\sur_k}$ for $1\leq k\leq n'$. 
\end{lemma}

\begin{proof}
Consider the decomposition in~\eqref{eq_Y decomposition}
\[
X_{(\underline a,\underline b)}=\pi^{-1}(Cone(L))\cup\pi^{-1}(P\setminus\{v_{n+2}\}). 
\]
Condition~\eqref{eqn_deg edge contraction condition} implies that
$\rho \colon P \to P'$ is a local homeomorphism near $v_{n+2}$. Hence, we may assume that $\rho(L)$ is a line segment that intersects $E'_{n'+1}$ and $E'_{n'+2}$ (see Figure~\ref{fig_line_segment}) and obtain a decomposition
\[
X_{(\underline a',\underline b')}=\pi'^{-1}(Cone(\rho(L)))\cup\pi'^{-1}(P'\setminus\{v'_{n'+2}\})
\]
for the orbit map $\pi'\colon X_{(\underline a',\underline b')}\to P'$. 
\begin{figure}
\begin{tikzpicture}
\node[fill=yellow!30, yellow!30, regular polygon, regular polygon sides=7, draw, minimum size = 2.5cm](m) at (0,0) {};

\node[opacity=0, regular polygon, regular polygon sides=7, draw, minimum size = 2.5cm](m) at (0,0) {};
\coordinate (1) at (m.corner 1); 
\coordinate (2) at (m.corner 2); 
\coordinate (3) at (m.corner 3); 
\coordinate (4) at (m.corner 4); 
\coordinate (5) at (m.corner 5); 
\coordinate (6) at (m.corner 6); 
\coordinate (7) at (m.corner 7); 

\draw[ yellow!30] (1)--(2)--(3)--(4)--(5)--(6)--(7)--cycle;

\draw (1)--(2)--(3);
\draw (4)--(5);
\draw (6)--(7);

\draw (1)--($0.2*(7)+0.8*(1)$);
\draw (7)--($0.8*(7)+0.2*(1)$);

\draw (3)--($0.2*(4)+0.8*(3)$);
\draw (4)--($0.8*(4)+0.2*(3)$);

\draw (5)--($0.2*(6)+0.8*(5)$);
\draw (6)--($0.8*(6)+0.2*(5)$);

\foreach \a in {1,...,7} { 
\draw[fill] (\a) circle (0.8pt); 
}

\draw[fill=red, red] ($0.5*(1)+0.5*(2)$) circle (1.5pt); 
\node[red, below right] at ($0.5*(1)+0.5*(2)$) {\scriptsize$L(1)$};
\draw[fill=red, red] ($0.5*(2)+0.5*(3)$) circle (1.5pt); 
\node[red, below right] at ($0.5*(2)+0.5*(3)$) {\scriptsize$L(0)$};
\draw[very thick, red] ($0.5*(1)+0.5*(2)$)--($0.5*(2)+0.5*(3)$);

\node[rotate=-28] at (m.side 7) {\scriptsize$\cdots$};
\node[rotate=130] at (m.side 3) {\scriptsize$\cdots$};
\node[rotate=51] at (m.side 5) {\scriptsize$\cdots$};

\node[left] at (m.side 2) {\scriptsize$E_{n+2}$};
\node[above left] at (m.side 1) {\scriptsize$E_{n+1}$};

\node[red] at ($0.6*(2)$) {$L$};
\node at ($0.1*(m.side 5)$) {$P$};

\draw[->] (2,0)--(4,0);
\node[above] at (3,0) {$\rho$};

\begin{scope}[xshift=180]
\node[fill=yellow!30, yellow!30, regular polygon, regular polygon sides=6, draw, minimum size = 2.5cm, rotate=30](n) at (0,0) {};

\node[opacity=0, regular polygon, regular polygon sides=6, draw, minimum size = 2.5cm, rotate=30](n) at (0,0) {};
\coordinate (1) at (n.corner 1); 
\coordinate (2) at (n.corner 2); 
\coordinate (3) at (n.corner 3); 
\coordinate (4) at (n.corner 4); 
\coordinate (5) at (n.corner 5); 
\coordinate (6) at (n.corner 6); 

\draw[ yellow!30] (1)--(2)--(3)--(4)--(5)--(6)--cycle;

\draw (1)--(2)--(3);
\draw (4)--(5)--(6);

\draw (1)--($0.2*(6)+0.8*(1)$);
\draw (6)--($0.8*(6)+0.2*(1)$);

\draw (3)--($0.2*(4)+0.8*(3)$);
\draw (4)--($0.8*(4)+0.2*(3)$);

\foreach \a in {1,...,6} { 
\draw[fill] (\a) circle (0.8pt); 
}

\draw[fill=red, red] ($0.5*(1)+0.5*(2)$) circle (1.5pt); 
\draw[fill=red, red] ($0.5*(2)+0.5*(3)$) circle (1.5pt); 
\draw[very thick, red] ($0.5*(1)+0.5*(2)$)--($0.5*(2)+0.5*(3)$);

\node[rotate=-30] at (n.side 6) {\scriptsize$\cdots$};
\node[rotate=150] at (n.side 3) {\scriptsize$\cdots$};

\node[left] at (n.side 2) {\scriptsize$E'_{n'+2}$};
\node[above left] at (n.side 1) {\scriptsize$E'_{n'+1}$};

\node[red, above] at ($0.2*(2)$) {\color{red}$\rho(L)$};
\node at ($0.2*(5)$) {$P'$};
\end{scope}
\end{tikzpicture}
\caption{A line segment in $P$ and an edge-contraction.}
\label{fig_line_segment}
\end{figure}

Consider the homotopy commutative diagram
\begin{equation}\label{diagram_edge contractions bwt skeleton 1}
\begin{tikzcd}
S^3\arrow[r,"\cong"]\arrow[d,"\phi"]
&
\pi^{-1}(L)\arrow[r,hook]
\arrow{d}{X(\rho,id_T)}
&\pi^{-1}(P\setminus\{v_{n+2}\})
\arrow{r}{\simeq}
\arrow{d}{X(\rho,id_T)}
&\bigvee^n_{\ell=1}S^2_\ell
\arrow{d}{p}\\
S^3\arrow[r,"\cong"]
&(\pi')^{-1}(\rho(L))\arrow[r,hook] 
&(\pi')^{-1}(P'\setminus\{v'_{n'+2}\})
\arrow{r}{\simeq}   
&\bigvee^{n'}_{k=1}S^2_{k}
\end{tikzcd}
\end{equation}
where $\phi$ is defined via the left square, the vertical maps in the middle square are restrictions of $X(\rho,id_T)$, and $p$ is the pinch map sending $S_{\sur_k}^2$ onto $S^2_k$ for $1\leq k\leq n'$ and contracting other $S^2_{\ell}$'s to the basepoint.
The left and middle squares commute trivially. The right square commutes up to homotopy since $\rho$ sends edges $E_{\sur_k}$ onto~$E'_{k}$ and contracts other edges to vertices.

Let $f\colon S^3\to\bigvee^n_{i=1}S^2_i$ and $f'\colon S^3\to\bigvee^{n'}_{k=1}S^2_k$ be the composites of maps in the top and the bottom rows in~\eqref{diagram_edge contractions bwt skeleton 1}, respectively. Then they are attaching maps of $4$-cells in $X_{(\underline a,\underline b)}$ and $X_{(\underline a',\underline b')}$. Hence \eqref{diagram_edge contractions bwt skeleton 1} can be extended to the homotopy commutative diagram
\begin{equation}\label{dgrm_homotopy cofibration}
\begin{tikzcd}
S^3\rar{f}\dar{\phi}   &\bigvee^n_{\ell=1}S^2_\ell
\rar{\jmath}\dar{p}    &X_{(\underline a,\underline b)}\dar{X(\rho,id_T)}\\
S^3\rar{f'}   
&\bigvee^{n'}_{k=1}S^2_{k}\rar{\jmath'}    &X_{(\underline a',\underline b')}
\end{tikzcd}
\end{equation}
where $\jmath$ is the composite $\bigvee^n_{i=1}S^2_i\simeq\pi^{-1}(P\setminus\{v_{n+2}\})\hookrightarrow X_{(\underline{a},\underline{b})}$ and $\jmath'$ is defined similarly. Since the rows are homotopy cofibration sequences, they induce a commutative diagram
\[
\begin{tikzcd}[column sep=scriptsize]
\cdots\rar{}    &H^{r-1}(S^3)\rar\dar{\phi^*}    &H^{r}(X_{(\underline a',\underline b')})\rar{(\jmath')^*}\dar{Y(\rho,id_T)^*}   
&\bigoplus^{n'}_{k=1}H^{r}(S^2_{k})
\rar{(f')^*}\dar{p^*}   &H^{r}(S^3)\rar{} \dar{\phi^*}   &\cdots\\
\cdots\rar{}    &H^{r-1}(S^3)\rar &H^{r}(X_{(\underline a,\underline b)})\rar{\jmath^*}    &\bigoplus^{n}_{\ell=1}H^{r}(S^2_\ell)\rar{f^*}   &H^{r}(S^3)\rar{} &\cdots
\end{tikzcd}
\]
where the rows are exact sequences. 

For $r=4$, generators $v\in H^4(X_{(\underline{a},\underline{b})})$ and $v'\in H^4(X_{(\underline{a'},\underline{b'})})$ correspond to~$Cone(L)$ and $Cone(\rho(L))$ respectively. Since $X(\rho,id)$ maps $Cone(L)$ homeomorphically onto~$Cone(\rho(L))$, we have $Y(\rho,id_T)^*(v')=\pm v$.  For $r=2$, each $(\jmath')^*(u'_k)$ 
is a generator of $H^2(S^2_k)$ for $1\leq k\leq n'$ and~$\jmath^\ast(u_{\ell})$ is a generator of $H^2(S^2_{\ell})$ for $1\leq \ell\leq n$. The commutativity of the middle square implies $X(\rho,id_T)^*(u'_k)=\pm u_{\sur_k}$.
\end{proof}

Edge contraction $\rho_i\colon P\to\triangle$ defined in~\eqref{eqn_rho_i} satisfies Condition~\eqref{eqn_deg edge contraction condition}, and hence Lemma~\ref{lemma_rho_S n cellular basis} applies to the cohomology morphism induced by the associated edge-contraction morphism.
We conclude this section by stating two useful lemmas that follow from this result.

\begin{lemma}\label{lemma_edge contraction gives cellular basis}
For $1\leq i\leq n$, let
$X(\rho_i,id)\colon X_{(\underline{a},\underline{b})}\to X_{(a_i,b_i)}$
be the edge contraction in~\eqref{eqn_rho_i}. Choose a generator $u^{\triangle}_i$ of $H^2(X_{(a_i,b_i)})$ for each $i$, and choose a generator $v$ of $H^4(X_{(\underline{a},\underline{b})})$. Then the set
\begin{equation}\label{eq_cell_basis_edge_cont}
\left\{ X(\rho_1,id_T)^*(u^{\triangle}_1),\ldots,X(\rho_n,id_T)^*(u^{\triangle}_n); ~ v\right\}
\end{equation}
forms a cellular basis for $\tilde{H}^*(X_{(\underline{a},\underline{b})})$.
\end{lemma}
\begin{proof}
Let $\{\tilde{u}^{\triangle}_i; \tilde{v}^{\triangle}_i\}$ and $\{\tilde{u}_1,\ldots,\tilde{u}_n; \tilde{v}\}$ be cellular bases for 
$\tilde{H}^*(X_{(a_i,b_i)})$ and $\tilde{H}^*(\xq{\underline a, \underline b})$, respectively. 
Since $H^2(X_{(a_i,b_i)})\cong\Z$, we have $u^{\triangle}_i=\pm\tilde{u}^{\triangle}_i$.
Hence, we apply Lemma~\ref{lemma_rho_S n cellular basis} to obtain
\[
X(\rho_i,id_T)^*(u^{\triangle}_i)=\pm\tilde{u}_i.
\]
Similarly, we have $v=\pm\tilde{v}$. 
Thus, the set \eqref{eq_cell_basis_edge_cont} forms a cellular basis for $\tilde{H}^*(X_{(\underline{a},\underline{b})})$.
\end{proof}

Given a primitive vector $(a,b)$, consider toric morphisms of degenerate toric spaces
\begin{equation}\label{eq_Xrho1_Xrho2}
X(\rho_1,id)\colon X_{(1,0),(a,b)}\to X_{(1,0)}
\quad\text{and}\quad
X(\rho_2,id)\colon X_{(a,b),(0,1)}\to X_{(0,1)}
\end{equation}
associated with the edge contraction $\rho_i\colon\square\to\triangle$ as in \eqref{eqn_rho_i}. Observe that $X_{(1,0),(a,b)}$ and  $X_{(a,b),(0,1)}$ are toric orbifolds if $b\neq 0$ and $a\neq 0$, respectively. In this case, we denote the corresponding cohomology isomorphisms in \eqref{eqn_wSR bar isom} by $\overline{\Phi}_1$ and $\overline{\Phi}_2$, respectively.

\begin{lemma}\label{lem_cel_basis_Y1001}
Let $\{y_1,\ldots,y_4\}$ be the canonical generators for $\SR(\square)$. There exist generators $u^{\triangle}_{(1,0)}\in H^2(X_{(1,0)})$ and $u^{\triangle}_{(0,1)}\in H^2(X_{(0,1)})$ such that
\begin{enumerate}
\item
if $X(\rho_1,id)$ is given in~\eqref{eq_Xrho1_Xrho2} and $b\neq0$, then $\overline{\Phi}_1\circ X(\rho_1,id)^\ast (u^{\triangle}_{(1,0)})=[by_2]$;
\item
if $X(\rho_2,id)$ is given in~\eqref{eq_Xrho1_Xrho2} and $a\neq0$, then $\overline{\Phi}_2\circ X(\rho_2,id)^\ast  (u^{\triangle}_{(0,1)})= [ay_1]$.
\end{enumerate}
\end{lemma}

\begin{proof}
Here, we prove the first claim in three steps. The second claim follows by a similar argument.

\medskip 
\noindent\underline{Step I:}\\
First, we assume $b=\pm 1$ and claim that any generator $\alpha\in H^2(X_{(1,0)})$ satisfies 
\begin{equation}\label{eq_alpha_lin_comb}
\left( \overline{\Phi}_1 \circ X(\rho_1,id)^*\right) (\alpha)=c[y_3]+d[y_4]
\end{equation}
with $d=\pm1$ and $c=0$ or $-2abd$.
Indeed, the assumption $b=\pm 1$ implies that~$X_{(1,0),(a,b)}$ is a smooth toric manifold. Hence, the isomorphism  $\overline{\Phi}_{1}$ maps $H^\ast(X_{(1,0),(a,b)})$ onto
\begin{equation}\label{eq_image_phi12}
\frac{\Z[y_1,y_2, y_3, y_4]}{\langle{y_1y_3,y_2y_4,y_1+ay_2+y_3,by_2+y_4}\rangle}\cong \frac{\mathbb{Z}[z_3, z_4]}{\left<y_4^2,ab y_3y_4-y_3^2\right>}.
\end{equation}
Therefore, $\left( \overline{\Phi}_1 \circ X(\rho_1,id)^*\right) (\alpha)$ can be written as a linear combination of $[y_3]$ and~$[y_4]$ as in \eqref{eq_alpha_lin_comb}. 

To calculate coefficients $c$ and $d$, we 
consider the edge-contraction 
\[
\rho^{\square}_{2}\colon\square\to\triangle
\]
shrinking $E_1$ of $\square$ to a point and let $\beta= \overline{\Phi}^{-1} ([ab y_1])$ be the generator of~$H^2(X_{(a,b)})$ for $\overline{\Phi} \colon H^\ast(X_{(a,b)}) \to \overline{\wSR}(\triangle, \lambda_{(a,b)})$ as in Example~\ref{eg_Xab}. 
Then Lemma~\ref{lemma_edge contraction gives cellular basis} implies that the set
\[
\{X(\rho_{1}, id_T)^\ast(\alpha), X(\rho^{\square}_{2}, id_T)^\ast(\beta)\}\] 
forms a basis of $H^2(X_{(1,0),(a,b)})$.
Moreover, it follows from~Lemma \ref{lemma_edge contract wSR} that 
\begin{equation*}
(\overline{\Phi}_{1} \circ X(\rho^{\square}_{2}, id_T)^\ast)(\beta)=\overline{\wSR}(\rho^{\square}_{2}, id_T)([ab y_1])=-[by_3].
\end{equation*}
Since $\overline{\Phi}_{1}$ is an isomorphism, the set $\{c[y_3]+d[y_4], -b[y_3]\}$ forms a basis of the degree~$2$ part of \eqref{eq_image_phi12}. Now, the assumption $b=\pm1$ implies that $d=\pm 1$.

Next, we recall from Example \ref{eg_X_(1,0)} that $X_{(1,0)}\simeq S^4\vee S^2$, so  $\alpha \cup \alpha=0$. 
Using~\eqref{eq_alpha_lin_comb} together with \eqref{eq_image_phi12}, we have  
\[
0=\alpha\cup\alpha=[(c^2ab+2cd)y_3y_4].
\]
Since $[y_3y_4]$ is a generator of $H^4(X_{(1,0),(a,b)})$ by Lemma~\ref{lem_5_deg_4_are_all_same}, we have $c=0$ or $cab=-2d$. 
If $a=0$ then $cab=-2d$ is impossible as we have shown $d=\pm1$. It must be that~$c=0$. If $a\neq0$ then $a=\pm1$ since $(a,b)$ is a primitive vector. Hence~$ab=\pm1$, and the equation $cab=-2d$ is equivalent to $c=-2abd$.

\medskip 

\noindent\underline{Step II:}\\
Here, we prove claim (1) for the case where $b=\pm1$.
Since the vector $(a,b)$ is primitive, we have $a=1,0$ or~$-1$. 
Let $P$ be a pentagon and $\lambda$ the characteristic function described below: 
\[
\begin{tikzpicture}
\node[regular polygon, regular polygon sides=5, draw, minimum size = 1.7cm](m) at (0,0) {};
\coordinate (1) at (m.corner 1); 
\coordinate (2) at (m.corner 2); 
\coordinate (3) at (m.corner 3); 
\coordinate (4) at (m.corner 4); 
\coordinate (5) at (m.corner 5); 

\node at ($1.4*(m.side 1)$) {\scriptsize$E_5$}; 
\node at ($1.4*(m.side 5)$) {\scriptsize$E_4$}; 
\node at ($1.3*(m.side 2)$) {\scriptsize$E_1$}; 
\node at ($1.3*(m.side 3)$) {\scriptsize$E_2$}; 
\node at ($1.3*(m.side 4)$) {\scriptsize$E_3$}; 

\draw[->] (2,0)--(3,0);
\node[above] at (2.5,0) {$\lambda$};

\begin{scope}[xshift=150]
\node[regular polygon, regular polygon sides=5, draw, minimum size = 1.7cm](m) at (0,0) {};
\coordinate (1) at (m.corner 1); 
\coordinate (2) at (m.corner 2); 
\coordinate (3) at (m.corner 3); 
\coordinate (4) at (m.corner 4); 
\coordinate (5) at (m.corner 5); 

\node at ($1.4*(m.side 1)$) {\footnotesize$(0,1)$}; 
\node at ($1.4*(m.side 5)$) {\footnotesize$(1,0)$}; 
\node at ($1.5*(m.side 2)$) {\footnotesize$(1,0)$}; 
\node at ($1.3*(m.side 3)$) {\footnotesize${(a',b)}$}; 
\node at ($1.5*(m.side 4)$) {\footnotesize$(0,b)$}; 
\end{scope}
\end{tikzpicture}
\]
where 
$a'=\begin{cases}
a	&\text{if }a\neq0\\
1	&\text{if }a=0.
\end{cases}$
Then $X(P,\lambda)=X_{(1,0),(a',b),(0,b)}$. Consider the commutative diagram
\begin{equation}\label{eq_pentagon_spaces}
\begin{tikzcd}[row sep=huge,column sep=large]
X_{(1,0),(a',b),(0,b)}\ar{r}{X(\rho_{12},id)}\ar{d}{X(\rho_{13},id)}	&X_{(1,0),(a',b)}\ar{d}{X(\rho_{12,1},id)}\\
X_{(1,0),(0,b)}\ar{r}{Y(\rho_{13,1},id)}								&X_{(1,0)}
\end{tikzcd}
\end{equation}
where $\rho_{ij}\colon P\to\square$ is the edge contraction given in~\eqref{eq_rho_ij}, and $\rho_{12,1}$ and $\rho_{13,1}$ are the edge contraction $\rho_1\colon\square\to\triangle$ shrinking the second edge of $\square$ to a point.
Diagram~\eqref{eq_pentagon_spaces}, together with the bottom face of \eqref{diagram_naturality of wSR n wSR bar}, induces the commutative diagram 
\begin{equation}\label{dgrm_pentagon induce cohmgly}
\small
\xymatrix{
	&H^\ast(X_{(1,0),(a',b),(0,b)})\ar@{-}[d]_-{\overline{\Phi}_{123}}	&	&H^*(X_{(1,0),(a',b)})\ar[ll]_-{X(\rho_{12},id)^*}\ar[dd]^-{\overline{\Phi}_{12}}\\
H^\ast(X_{(1,0),(0,b)})\ar[dd]_-{\overline{\Phi}_{13}}\ar[ur]^-{X(\rho_{13},id)^*}	&\ar[d]	&H^*(X_{(1,0)})\ar[ll]^(0.3){X(\rho_{13,1},id)^*}\ar[ur]_-{X(\rho_{12,1},id)^*}	&\\
	&\displaystyle\frac{\Z[y_1,y_2,y_3,y_4,y_5]}{\mathcal{I}_{123}}	&	&\displaystyle\frac{\Z[y_1,y_2,y_3,y_4]}{\mathcal{I}_{12}}\ar[ll]^-{\overline{\wSR}(\rho_{12},id)}\\
\displaystyle\frac{\Z[y_1,y_2,y_3,y_4]}{\mathcal{I}_{13}}\ar[ur]_-{\overline{\wSR}(\rho_{13},id)}&	&	&
}
\end{equation}
where $\overline{\Phi}_{12},\overline{\Phi}_{13}$ and $\overline{\Phi}_{123}$ are isomorphisms in~\eqref{eqn_wSR bar isom}, and $\mathcal{I}_{12},\mathcal{I}_{13}$ and $\mathcal{I}_{123}$ are ideals
\begin{align*}
\mathcal{I}_{12}&=\langle{y_1y_3,y_2y_4,y_1+a'y_2+y_3,by_2+y_4}\rangle,\\
\mathcal{I}_{13}&=\langle{y_1y_3,y_2y_4,y_1+y_3,by_2+y_4}\rangle,\\
\mathcal{I}_{123}&=\langle{y_1y_3, y_1y_4, y_2y_4, y_2y_5, y_3y_5, y_1+a'y_2+y_4, by_2+by_3+y_5}\rangle.
\end{align*}
Take a generator $\alpha \in H^2(X_{(1,0)})$.
The commutativity of~\eqref{dgrm_pentagon induce cohmgly} implies
\begin{equation}\label{eq_parallelogram_equation}
\small
(\overline{\wSR}(\rho_{12}, id)\circ \overline{\Phi}_{12} \circ X(\rho_{12,1}, id)^\ast)(\alpha)=
(\overline{\wSR}(\rho_{13}, id)\circ \overline{\Phi}_{13} \circ X(\rho_{13,1}, id)^\ast)(\alpha).
\end{equation}
From the Step I discussion, we know that for some $d,d'\in\{1,-1\}$ and $\epsilon\in\{0,1\}$
\[\begin{cases}
\left(\overline{\Phi}_{12}\circ X(\rho_{12,1},id)^*\right) (\alpha)=d[y_4]-2a'bd\epsilon[y_3],\\
\left( \overline{\Phi}_{13} \circ X(\rho_{13,1},id)^*\right) (\alpha)=d'[y_4].
\end{cases}
\]
Substitute them into~\eqref{eq_parallelogram_equation} and apply Lemma~\ref{lemma_rho_ij to sq} to obtain
\[
d'[y_5]=d[y_5]-2a'bd\epsilon([y_4]-a'[y_3]).
\]
If $\epsilon=1$ then it follows that $H^2(X_{(1,0),(a',b),(b,0)})$ has two generators which is false.
Therefore $\epsilon=0$ and $d=d'$.

Let $u^{\triangle}_{(1,0)}=-d\alpha$. Then it is a generator of $H^2(X_{(1,0)})$ as $d=\pm1$. Due to the relation $[by_2+y_4]=0$ in $\Z[y_1,\ldots,y_4] / \mathcal{I}_{12}$ and $\Z[y_1,\ldots,y_4] / \mathcal{I}_{13}$, 
we have
\[
\begin{cases}
\left(\overline{\Phi}_{12}\circ  X(\rho_{12,1},id)^*\right) (u^{\triangle}_{(1,0)})=[-y_4]=[by_2]\\
\left( \overline{\Phi}_{13}\circ X(\rho_{13,1},id)^*\right) (u^{\triangle}_{(1,0)})=[-y_4]=[by_2].
\end{cases}
\]
Hence the lemma is proved for the case where $b=\pm1$.

\medskip 
\noindent\underline{Step III:}\\
We finally prove the lemma in the general case.
It suffices to show that the generator~$u^{\triangle}_{(1,0)}$ constructed in Step II satisfies
\begin{equation}\label{eqn_rho_1 image of u_(1,0)}
X(\rho_1,id)^*(u^{\triangle}_{(1,0)})=\overline{\Phi}^{-1}_1([by_2])\in H^2(X_{(1,0),(a,b)})
\end{equation}
for any primitive vector $(a,b)\in\Z^2$ with $b\neq0$. Due to Step II, we assume $b\neq\pm1$. Then $a\neq0$ as $(a,b)$ is primitive. 

Let $\epsilon_1$ and $\epsilon_2$ be the signs of $a$ and $b$, and let 
\[
X(id,\sigma)\colon X_{(1,0),(a,b)}\to X_{(1,0),(\epsilon_1,\epsilon_2)}\] 
be the rescaling morphism in Definition~\ref{def_rescaling}
with respect to $(a,b)$, where $\sigma \colon T^2\to T^2$ is given by 
$\sigma(t_1, t_2)= (t_1^{|b|}, t_2^{|a|})$. 
Then, there is a commutative diagram
\[
\xymatrix{
H^*(X_{(1,0)})\ar[rr]^-{X(\rho_1,id)^*}\ar[d]_-{X(id,\sigma)^*}	&&H^*(X_{(1,0),(\epsilon_1,\epsilon_2)})\ar[d]^-{X(id,\sigma)^*}\\
H^*(X_{(1,0)})\ar[rr]^-{X(\rho_1,id)^*}								&&H^*(X_{(1,0),(a,b)})}
\]
where the horizontal maps are induced by the edge contraction $X(\rho_1,id)$ in~\eqref{eqn_rho_i}, and the left vertical map $X(id,\sigma)^*$ is the restriction of the right $X(id,\sigma)^*$.  

Consider the images of $u^{\triangle}_{(1,0)}$ under the composites of maps round the top right and the bottom left corners.
On one hand, by Lemmas~\ref{lemma_rescaling morphism wSR} and~\ref{lemma_rho_ij to sq} and the Step~I discussion,  we have
\[
\left(X(id,\sigma)^*\circ X(\rho_i,id)^*\right)(u^{\triangle}_{(1,0)})=
\left(X(id,\sigma)^*\circ\overline{\Phi}^{-1}\right)([\epsilon_2z_2])=\overline{\Phi}^{-1}([|a|bz_2]).
\]
On the other hand, by Lemma~\ref{lemma_rescale Y_(1,0)} we have
\[
\left(X(\rho_i,id)^*\circ X(id,\sigma)^*\right)(u^{\triangle}_{(1,0)})=X(\rho_i,id)^*(|a|u^{\triangle}_{(1,0)}).
\]
Since $a\neq0$, the commutativity of the diagram implies~\eqref{eqn_rho_1 image of u_(1,0)}. Therefore the proof is complete.
\end{proof}

\section{Proof of the main theorem}
\label{sec_pf_main_thm}

Let $X(P,\lambda)$ be a $4$-dimensional toric orbifold with a smooth vertex.
Below, we define a cellular basis of  $H^*(X(P,\lambda))$, called the \emph{algebraic cellular basis}, and prove that the cup products of its elements satisfy Equation~\eqref{eqn_cup prod equation}.

Relabeling vertices of $P$ and changing basis in $\Z^2$, we assume $X(P,\lambda)=X_{(\underline{a},\underline{b})}$ as in~\eqref{eq_1001_char_fun}.
Recall from \eqref{eqn_wSR bar isom} that there is a ring isomorphism
\[
\overline{\Phi}_X\colon H^*(X(P,\lambda))\to\overline{\wSR}(P,\lambda).
\]

\begin{dfn}\label{dfn_good cellular basis}
Let $\{y_1,\ldots,y_{n+2}\}$ be the canonical generators of $\SR(P)$. 
The \emph{algebraic cellular basis} for $\tilde{H}^*(X_{(\underline{a},\underline{b})})$ is a set $\mathcal{B}=\{ u_1,\ldots,u_n,v\}$ of cohomology classes defined by 
\begin{align*}
    u_i&=\overline{\Phi}_X^{-1}\left(\left[
\sum^{i-1}_{k=1}a_kb_iy_k+a_ib_iy_i+\sum^n_{l=i+1}a_ib_ly_l\right]\right)\in H^2(X_{(\underline{a},\underline{b})});\\
v&=\overline{\Phi}^{-1}_X([y_{n+1}y_{n+2}])\in H^4(X_{(\underline{a},\underline{b})}).
\end{align*}
\end{dfn}

Before showing that $\mathcal{B}$ is a cellular basis and proving Theorem~\ref{thm_main}, we verify that $u_i$ and $v$ are well-defined. By Lemma~\ref{lem_5_deg_4_are_all_same}, $v$ is a generator of $H^4(X_{\underline{a},\underline{b}})$. The well-definedness of $u_i$ is subtler, since it is not a priori clear that the polynomials defining $u_i$ lie in $\wSR(P,\lambda)$. 
The following lemma shows that $u_i$ arises as the image of a cohomology class in $H^2(X_{(a_i,b_i)})$ under an edge contraction, and is therefore well-defined.

\begin{lemma}\label{lemma_describe u_i}
For $1\leq i\leq n$, let $X(\rho_i,id)\colon X_{(\underline{a},\underline{b})}\to X_{(a_i,b_i)}$ be the edge contraction morphism in~\eqref{eqn_rho_i}.
Then the cohomology class $u_i$ in Definition~\ref{dfn_good cellular basis} satisfies
\[
u_i=X(\rho_i,id)^*(u^{\triangle}_{(a_i,b_i)}),
\]
where $u^{\triangle}_{(a_i,b_i)}\in H^2(X_{(a_i,b_i)})$ is a generator given as follows:
\begin{itemize}
\item
if $a_ib_i\neq0$ then
$u^{\triangle}_{(a_i,b_i)}=\overline{\Phi}^{-1}([a_ib_iy_1])$ is the generator as in Example~\ref{eg_Xab};
\item
$u^{\triangle}_{(1,0)}$ and $u^{\triangle}_{(0,1)}$ are the generators given in Lemma~\ref{lem_cel_basis_Y1001};
\item
$u^{\triangle}_{(-1,0)}=-u^{\triangle}_{(1,0)}$ and $u^{\triangle}_{(0,-1)}=-u^{\triangle}_{(0,1)}$.
\end{itemize}
\end{lemma}

\begin{proof}
It suffices to show that $\overline{\Phi}\circ X(\rho_i,id)^*(u^{\triangle}_{(a_i,b_i)})$ and $\overline{\Phi}(u_i)$ have the same images in $\overline{\wSR}[P,\lambda]$. The proof is divided into  4 cases.\\

\noindent\underline{Case 1}: $(a_i,b_i)$ satisfies $a_ib_i\neq0$.\\
Use the bottom commutative square in~\eqref{diagram_naturality of wSR n wSR bar} and Lemma~\ref{lemma_edge contract wSR} to obtain
\begin{align*}
\overline{\Phi}\circ X(\rho_i,id)^*(u^{\triangle}_{(a_i,b_i)})
&=\overline{\wSR}(\rho_i,id)^*([a_ib_ix_1])\\
&=\left[\sum^{i-1}_{k=1}a_kb_iy_k+a_ib_iy_i+\sum^n_{k=i+1}a_ib_ky_k\right]\\
&=\overline{\Phi}(u_i).
\end{align*}

\noindent\underline{Case 2}: $(a_i,b_i)=(1,0)$.\\
Since $X_{(\underline{a},\underline{b})}$ is a toric orbifold,  the characteristic vector $\lambda(E_{i+1})=(a_{i+1},b_{i+1})$ satisfies $b_{i+1}\neq0$.
First, we assume $b_{i+1}=\pm1$.
Then $X_{(1,0),(a_{i+1},b_{i+1})}$ is a toric orbifold.
The edge contraction $\rho_i\colon P\to \triangle$ factors as a composite of edge contractions
\begin{equation*}\label{eqn_rho_ij,k}
\rho_i\colon P\xrightarrow{\rho_{i,i+1}}\square\xrightarrow{\tilde{\rho}_1}\triangle.
\end{equation*}
where $\rho_{i,i+1}$ is given in~\eqref{eq_rho_ij} and $\tilde{\rho}_1$ shrinks the second edge of $\square$ to a point.
Apply
Lemmas~\ref{lemma_rho_ij to sq} and~\ref{lem_cel_basis_Y1001} to obtain
\begin{align*}
\left(\overline{\Phi}\circ X(\rho_i,id)^*\right)(u^{\triangle}_{(1,0)})
&=\left(\overline{\Phi}\circ X(\rho_{i,i+1},id)^*\circ X(\tilde{\rho}_1,id)^*\right)(u^{\triangle}_{(1,0)})\\
& \left(\overline{\Phi}\circ X(\rho_{i,i+1},id)^*\circ\overline{\Phi}^{-1}_{\square}\right)([b_{i+1}y_2])\\
&=\overline{\wSR}(\rho_{i, i+1}, id)([b_{i+1} y_2]) \\
&=[b_{i+1}y_{i+1}+\cdots+b_ny_n]\\
&=\overline{\Phi}(u_i)
\end{align*}
where $\overline{\Phi}_{\square}$ is the isomorphism~\eqref{eqn_wSR bar isom} for $H^*(X_{(1,0),(a_{i+1},b_{i+1})})$.
Therefore, the lemma holds when $(a_i,b_i)=(1,0)$ and $b_{i+1}=\pm1$.

Next, assume $|b_{i+1}|>1$. Since $(a_{i+1},b_{i+1})$ is a primitive vector, $a_{i+1}$ is non-zero. Consider the rescaling morphism
\[
X(id,\sigma)\colon X(P,\lambda)\to X(P,\lambda')
\]
induced by $\sigma(t_1,t_2)=(t_1^{|b_{i+1}|},t_2^{|a_{i+1}|})$, where
\[
\lambda'(E_j)=\left(\frac{|b_{i+1}|a_j}{g_{j,i+1}},\frac{|a_{i+1}|b_j}{g_{j,i+1}}\right)
\quad\text{and}\quad
g_{j,i+1}=\gcd(|a_jb_{i+1}|,|b_ja_{i+1}|)
\]
Denote $\lambda'(E_j)=(a'_j,b'_j)$
and let $(\underline{a}',\underline{b}')=\{(a'_j,b'_j)\}^n_{j=1}$. There is a commutative diagram
\[
\xymatrix{
H^*(X_{(1,0)})\ar[rr]^-{X(\rho_i,id)^*}\ar[d]^-{X(id,\sigma)^*}	&&H^*(X_{(\underline{a}',\underline{b}')})\ar[d]^-{X(id,\sigma)^*}\\
H^*(X_{(1,0)})\ar[rr]^-{X(\rho_i,id)^*}					&&H^*(X_{(\underline{a},\underline{b})}).
}
\]
Since $b'_{i+1}=\pm1$, the above result together with Lemmas~\ref{lemma_rescaling morphism wSR} and~\ref{lemma_rescale Y_(1,0)} gives us the following: 
\begin{align*}
\left( \overline{\Phi}\circ X(\rho_i,id)^*\circ X(id,\sigma)^*\right) (u^{\triangle}_{(1,0)})
&=\left( \overline{\Phi}\circ X(id,\sigma)^*\circ X(\rho_i,id)^*\right) (u^{\triangle}_{(1,0)}); \\
\left( \overline{\Phi}\circ X(\rho_i,id)^*\right) (|a_{i+1}|u^{\triangle}_{(1,0)})
&=[g_{i+1,i+1}b'_{i+1}y_{i+1}+\cdots+g_{i+1,n}b'_ny_n]\\
&=[|a_{i+1}|b_{i+1}y_{i+1}+\cdots+|a_{i+1}|b_ny_n]\\
&=|a_{i+1}|\overline{\Phi}(u_i).
\end{align*}
Since $a_{i+1}\neq0$, we can divide both sides by $|a_{i+1}|$, which establishes the claim in this case.\\

\noindent\underline{Case 3}: $(a_i,b_i)=(-1,0)$.

Let $\lambda''$ be the characteristic function on $P$ given by
\[
\lambda''(E_j)=\begin{cases}
\lambda(E_j)	&\text{if }j\neq i\\
(1,0)		&\text{if }j= i,
\end{cases}
\]
and let $\overline{\Phi}''\colon H^*(X(P,\lambda''))\to\overline{\wSR}(P,\lambda'')$ the corresponding isomorphism in~\eqref{eqn_wSR bar isom}. From the discussion of Case~2 and Lemma~\ref{lem wSR isom for P with -lambda} we have
\begin{align*}
\left( \overline{\Phi}\circ X(\rho_i,id)^*\right) (u^{\triangle}_{(-1,0)})
&=\left( \overline{\wSR}(id,id)\circ\overline{\Phi}''\circ X(\rho_i,id)^*\right) (-u^{\triangle}_{(1,0)})\\
&=\overline{\wSR}(id,id)([-b_{i+1}y_{i+1}-\cdots-b_ny_n])\\
&=[-b_{i+1}y_{i+1}-\cdots-b_ny_n]\\
&=\overline{\Phi}(u_i).
\end{align*}

\noindent\underline{Case 4}: $(a_i,b_i)=(0,\pm1)$\\
The case can be proved using arguments similar to those in Cases 2 and 3.
\end{proof}

Next, we show that $\mathcal{B}$ is a cellular basis for $\tilde{H}^*(X_{(\underline{a},\underline{b})})$, which justifies its name.

\begin{lemma}\label{lemma_alg cellular basis is cellular basis}
The  algebraic cellular basis $\mathcal{B}=\{u_1,\ldots,u_n,v\}$ given in Definition~\ref{dfn_good cellular basis} is a cellular basis for $\tilde{H}^*(X_{(\underline{a},\underline{b})})$.
\end{lemma}

\begin{proof}
The lemma follows immediately from Lemmas~\ref{lem_5_deg_4_are_all_same}, \ref{lemma_edge contraction gives cellular basis}, and~\ref{lemma_describe u_i}.
\end{proof}

Now we give the proof of the main theorem.

\begin{proof}[Proof of Theorem~\ref{thm_main}]
By Lemma~\ref{lemma_alg cellular basis is cellular basis} the algebraic cellular basis $\mathcal{B}=\{u_1,\ldots,u_n,v\}$ spans $\tilde{H}^*(X_{(\underline{a},\underline{b})})$. We claim that $u_i\cup u_j=a_ib_jv.$ for $1\leq i\leq j\leq n$.

For convenience write $u_i=\overline{\Phi}^{-1}_X(\alpha_i)$ where
\[
\alpha_i=\left[\sum^{i-1}_{k=1} a_kb_iy_k+a_ib_iy_i+\sum^n_{l=i+1}a_ib_ly_l\right]\in\overline{\wSR}(P,\lambda).
\]
Since $\overline{\Phi}_X$
is an isomorphism, it suffices to show that for $1\leq i\leq j\leq n$ 
\begin{equation}\label{eqn_alg pf of main thm}
\alpha_i\cdot\alpha_j=[a_ib_jy_{n+1}y_{n+2}]
\end{equation}
in $\overline{\wSR}[P,\lambda]$.

First, using the relation $[b_1y_1+\cdots+b_ny_n+y_{n+2}]=0$ we have
\begin{align*}
\alpha_i
&=\left[b_i\sum^{i-1}_{k=1}a_ky_k+a_ib_iy_i+a_i\sum^n_{l=i+1}b_ly_l\right]\\
&=\left[b_i\sum^{i-1}_{k=1}a_ky_k+a_ib_iy_i-a_i\left(\sum^i_{l=1}b_ly_l+y_{n+2}\right)\right]\\
&=\left[\sum^{i-1}_{k=1}c_ky_k-a_iy_{n+2}\right]
\end{align*}
where we write $c_k=a_kb_i-a_ib_k$.
Similarly, the relation $[a_1y_1+\cdots+a_ny_n+y_{n+1}]=0$ implies
\[
\alpha_j=\left[\sum^{n}_{k=j+1}d_ky_k-b_jy_{n+1}\right]
\]
where $d_k=a_kb_j-a_jb_k$.
Their product equals
\begin{align*}
\alpha_i\cdot\alpha_j&=\left[a_ib_jy_{n+1}y_{n+2}+\left(\sum^{i-1}_{k=1}c_ky_k\right)\cdot\left(\sum^{n}_{l=j+1}d_ly_l\right)\right.\\
&\qquad\left.-\left(\sum^{i-1}_{k=1}c_ky_k\right)b_jy_{n+1}-\left(\sum^{n}_{l=j+1}d_ly_l\right)a_iy_{n+2}\right]\\
&=[a_ib_jy_{n+1}y_{n+2}].
\end{align*}
In the first line, the product $\left(\sum^{i-1}_{k=1}c_ky_k\right)\cdot\left(\sum^{n}_{l=j+1}d_ly_l\right)$ is a polynomial in
\[
\{y_ky_l\mid
1\leq k\leq i-1, j+1\leq l\leq n\}.
\]
Since $i<j$, the indices $k,l$ satisfy $|k-l|>1$, implying that each $[y_ky_l]=0$ in~$\overline{\wSR}[P,\lambda]$. Hence the product is zero.
A similar argument shows that the products in the second lines are zero.
Therefore, all terms on the right hand side of the first equality are zero except $a_ib_jy_{n+1}y_{n+2}$. Equation~\eqref{eqn_alg pf of main thm} holds and the theorem follows.
\end{proof}

\begin{eg}
Let $P$ be a square and let $\lambda$ be the characteristic function as shown in Figure~\ref{fig_eg main thm}.
Since $v_4$ is a smooth point and $\lambda$ satisfies~\eqref{eqn_lambda refined form},
Theorem~\ref{thm_main} implies that $\tilde{H}^*(X(P,\lambda))$ has the algebraic cellular basis
\[
\mathcal{B}=\{u_1,u_2;v\}
\]
where the degree-$2$ generators $u_i$ and degree-$4$ generator $v$ are
\begin{itemize}
\item
$u_1=\overline{\Phi}^{-1}([-2x_1-2x_2])$ and $u_2=\overline{\Phi}^{-1}([4x_1-2x_2])$;
\item
$v=\overline{\Phi}^{-1}([x_3x_4])$.
\end{itemize}
Following Theorem \ref{thm_main}, the cellular cup product representation (Definition \ref{def_cellular_basis_cel_cup_repn}) with respect to $\mathcal{B}$ is
\[
M(\mathcal{B},X(P,\lambda))=\begin{pmatrix}
-2	&4\\
4	&-2
\end{pmatrix}.
\]

\begin{figure}
\begin{tikzpicture}[baseline=(current bounding box.center)]
\node[regular polygon, regular polygon sides=4, draw, minimum size = 1.5cm](Q) at (0,0) {};
\node[above] at (Q.side 1) {\scriptsize$E_3$}; 
\node[left] at (Q.side 2) {\scriptsize$E_4$}; 
\node[below] at (Q.side 3) {\scriptsize$E_1$}; 
\node[right] at (Q.side 4) {\scriptsize$E_2$}; 

\draw[->] (1.8,0)--(2.8,0);
\node[above] at (2.3,0) {$\lambda$};

\node[regular polygon,  regular polygon sides=4, draw, minimum size = 1.5cm](Qi) at (5,0) {};
\node[above] at (Qi.side 1) {\scriptsize$(1,0)$}; 
\node[left] at (Qi.side 2) {\scriptsize$(0,1)$}; 
\node[below] at (Qi.side 3) {\scriptsize$(-2,1)$}; 
\node[right] at (Qi.side 4) {\scriptsize$(1,-2)$}; 
\end{tikzpicture}
\caption{Characteristic function on a square.}
\label{fig_eg main thm}
\end{figure}
\end{eg}

\appendix
\section{Orientation of toric orbifolds}
Using the quotient construction 
$X(P, \lambda)=P\times {T^d}/_\sim$
from Definition~\ref{dfn_toric_orb}, we define an orientation on $X(P, \lambda)$ as follows. Equip $P$ and $T^d$ with their standard orientations as subspaces of $\mathbb{R}^d$ and $\mathbb{C}^d$, respectively, and choose the orientation on~$X(P, \lambda)$ such that the quotient map 
\[
\pi\colon P\times T^d\to X(P,\lambda)
\] 
reverses orientation if $d=4k+1$ or $4k+2$, and preserves orientation if $d=4k$ or $4k+3$. We call it the \emph{cellular orientation} of $X(P, \lambda)$. When $X(P,\lambda)$ is a toric variety, this orientation agrees with that induced by the canonical orientation of its open dense subset $(\C^*)^d\subset X(P,\lambda)$.

To be more precise, given a $d$-dimensional simple rational polytope $P$ in $M\otimes_{\Z}\R$ where $M$ is an integral lattice, one
can define a projective toric variety $X_P$ of complex dimension $d$ (see~\cite[Chapter 2]{CLS}).
Due to the work of~\cite{Jur} there is a homeomorphism 
\begin{equation}\label{eq_underlying_homeo}
h\colon X_P \to X(P, \lambda) =P\times T^d /_\sim,
\end{equation}
where $\lambda$ is the characteristic function that assigns to each facet of $P$ the primitive vector orthogonal to its supporting hyperplane. See for instance~\cite[Section 1]{Fra} and the references therein. 
As a toric variety, $X_P$ carries an orientation induced by the natural orientation on its open dense subset $(\C^\ast)^d$.
Moreover, $h$ restricts to a homeomorphism from $(\C^\ast)^d\subset X_P$ onto $\interior(P) \times T^d$ where $\interior(P)$ denotes the interior of $P$. Under the identification $(\C^*)^d\cong\R_{>0}^d\times T^d$ given by
\[
(r_1e^{i\theta_1},\ldots, r_de^{i \theta_d})\mapsto(r_1,\ldots, r_d, \theta_1,\ldots, \theta_d),
\]
this map is orientation-reversing if $d\equiv 1$ or $2\pmod{4}$ and is orientation-preserving if $d\equiv0$ or $3\pmod{4}$. Hence, the canonical orientation of $X_P$ as a toric variety agrees with the cellular orientation of $X(P,\lambda)$ defined above.

Recall that a $4$-dimensional toric orbifold $X(P, \lambda)$ satisfies $H_4(X(P, \lambda))\cong \Z$. An orientation of $X(P, \lambda)$ determines a generator of $H_4(X(P, \lambda))$, and vise versa. Such a generator is called the \emph{fundamental class}, which we denote by $[X(P, \lambda)]$. The following proposition shows that the algebraic cellular basis $v\in H^4(X(P, \lambda))$ in Definition \ref{dfn_good cellular basis} is the Kronecker dual of the fundamental class with respect to the cellular orientation on $X(P, \lambda)$.

\begin{prop}\label{prop_fundamental class dual}
Let $X(P,\lambda)$ be a $4$-dimensional toric orbifold equipped with the cellular orientation, and let $[X(P,\lambda)]$ be its fundamental class. If $X(P,\lambda)$ has a smooth vertex 
and $v\in H^4(X(P,\lambda))$ is the degree-$4$ generator given in Definition~\ref{dfn_good cellular basis}, then
\[
\langle v,[X(P,\lambda)]\rangle=1,
\]
where $\langle~,~\rangle\colon H^i(X(P,\lambda))\otimes H_i(X(P,\lambda))\to\Q$ is the natural pairing.
\end{prop}

\begin{proof}
Here we work with rational coefficients and suppress them from the cohomology notation for brevity.

First, consider the special case of $X_{(-1,-1)}\cong\C\PP^2$.
As a toric variety, its underlying fan $\Sigma$
has $1$-dimensional cones $\Sigma^{(1)}=\{\rho_1,\rho_2,\rho_3\}$ where
\[
\rho_1=\{(-t,-t)\mid t\geq0\},\quad
\rho_2=\{(t,0)\mid t\geq0\},\quad
\rho_3=\{(0,t)\mid t\geq0\}.
\]
Let $D_{\rho_i}$ be the Weil divisor of $\rho_i$, and let $v^{\triangle}\in H^4(X_{(-1,-1)})$ be the degree-$4$ generator given in Definition~\ref{dfn_good cellular basis}. Then we have
\[
\overline{\Phi}(v^{\triangle})=[y_2y_3]=\PD([D_{\rho_2}])\cdot \PD([D_{\rho_3}])=\PD([pt])
\]
where $\overline{\Phi}$ is the isomorphism~\eqref{eqn_wSR bar isom}, $\PD\colon H_i(X_{(-1,-1)})\to H^{4-i}(X_{(-1,-1)})$ is the Poincar\'{e} duality map, 
and $[pt]\in H_0(X_{(-1,-1)})$ is the homology class representing a point in $X_{(-1,-1)}$.
Then we have
\[
\langle v^{\triangle},X_{(-1,-1)}]\rangle=
[X_{(-1,-1)}]\cap v^{\triangle}
=\PD^{-1}(\PD([pt]))
=[pt]=1.
\]

Second, we consider the case of $X_{(1,1)}$. Since  $X(id,id)\colon X_{(1,1)}\to X_{(-1,-1)}$ is the identity map, Diagram~\eqref{dgrm_wSR isom for P with -lambda diagram} implies that $v^{\triangle}$ is also the degree-$4$ generator of the algebraic cellular basis for $H^4(X_{(1,1)})$. Moreover, $[X_{(1,1)}]$ equals $[X_{(-1,-1)}]$ as the cellular orientations on $X_{(1,1)}$ and $X_{(-1,-1)}$ coincide. Therefore, the claim also holds for $X_{(1,1)}$.

Next, we prove for the general case. Let $X(P,\lambda)$ be a $4$-dimensional toric orbifold having a smooth vertex, which we take to be $v_{n+2}$ without loss of generality. Suppose that there exists an $i\in\{1,\ldots,n\}$ such that $\lambda(E_i)=\epsilon(1,1)$ for $\epsilon=\pm1$. Consider the edge-contraction morphism
\[
X(\rho_i,id)\colon X(P,\lambda)\to X_{\epsilon(1,1)}
\]
where $\rho_i$ is given in~\eqref{eqn_rho_i}.
Let $v\in H^4(X(P,\lambda))$ be the degree-$4$ generator given in Definition~\ref{dfn_good cellular basis}.
By construction and Lemma~\ref{lem_5_deg_4_are_all_same} we have $X(\rho_i,id)^*(v^{\triangle})=v$. Since $X(\rho_i,id)$ preserves orientation, we have $X(\rho_i,id)_*([X(P,\lambda)])=[X_{\epsilon(1,1)}]$ and
\[
\langle v,[X(P,\lambda)]\rangle
=\langle X(\rho_i,id)^*(v^{\triangle}),X(\rho_i,id)_*^{-1}([X_{\epsilon(1,1)}])\rangle
=\langle v^{\triangle},[X_{\epsilon(1,1)}]\rangle
=1.
\]
Otherwise all characteristic vectors $\lambda(E_i)\neq\pm(1,1)$. Let $P'$ be an $(n+3)$-gon and let $\lambda'$ be a characteristic function on $P'$ given by
\[
\lambda'(E_i)=\begin{cases}
\lambda(E_i)		&\text{for }1\leq i\leq n,\\
(1,1)			&\text{for }i=n+1,\\
(1,0)			&\text{for }i=n+2,\\
(0,1)			&\text{for }i=n+3.
\end{cases}
\]
Let $v'\in H^4(X(P',\lambda'))$ be the degree-$4$ generator given in Definition~\ref{dfn_good cellular basis}. From the above discussion $\langle v',[X(P',\lambda')]\rangle=1$. 
Consider the edge-contraction morphism
\[
X(\rho',id)\colon X(P',\lambda')\to X(P,\lambda)
\]
where $\rho'\colon P'\to P$ shrinks $E'_{n+1}$ to a point. Then $X(\rho',id)^*(v)=v'$ and
\[
\langle v,[X(P,\lambda)]\rangle=\langle X(\rho',id)^*(v),X(\rho',id)^{-1}_*([X(P,\lambda)])\rangle=
\langle v',[X(P',\lambda')]\rangle=1,
\]
where the second equality is due to the fact that $X(\rho',id)$ preserves orientation.
\end{proof}

\section{Partition and connected sums}

Here we discuss an operation on toric orbifolds induced by partitioning a polygon into two. Although it is not a toric morphism in the sense of Definition~\ref{def_toric_morphism}, the underlying idea is similar. Both involve combinatorial operations on polygons to produce new toric orbifolds and to obtain information about the original toric orbifold.

\begin{dfn}\label{dfn_connected sum}
Let $(P,\lambda)$ be a characteristic pair where $P$ is an $(n+2)$-gon for $n\geq 2$. A \emph{partition} of $(P,\lambda)$ is a pair $\{(P^+,\lambda^+),(P^-,\lambda^-)\}$ of two degenerate characteristic pairs such that
\begin{itemize}
\item $P^+$ is an $(r+2)$-gon and $P^-$ is an $(s+2)$-gon with $r+s=n$;
\item writing $\{E_i^+\mid 1\leq i \leq r+2\}$ and  $\{E_i^- \mid 1\leq i \leq s+2\}$ for facets of $P^+$ and~$P^-$, respectively, $\lambda^+$ and $\lambda^-$ are defined by 
\begin{align*}
\lambda^+(E_i^+)&= 
\begin{cases} 
\lambda(E_i) & i=1, \dots, r+1;\\
\lambda(E_{n+2}) & i=r+2,
\end{cases}\\
\lambda^-(E_j^-)&= 
\begin{cases} 
\lambda(E_{r+1+j})& j=1, \dots, s+1;\\
\lambda(E_{r+1}) & j=s+2.
\end{cases}
\end{align*}
\end{itemize}
\end{dfn}

Intuitively, $P^+$ and $P^-$ are obtained by cutting $P$ along a line from $E_{n+2}$ and~$E_{r+1}$ and attaching new vertices to the resulting half-edges to close the boundaries. See Figure~\ref{fig_partition_polygon}.

\begin{figure}
\begin{tikzpicture}
\coordinate (A) at (22.5:1.3);
\coordinate (B) at (67.5:1.3);
\coordinate (C) at (112.5:1.3);
\coordinate (D) at (157.5:1.3);
\coordinate (E) at (202.5:1.3);
\coordinate (F) at (247.5:1.3);
\coordinate (G) at (292.5:1.3);
\coordinate (H) at (337.5:1.3);

\draw[thick, fill=yellow!30] (A)--(B)--(C)--(D)--(E)--(F)--(G)--(H)--cycle;

\foreach \P in {A,B,C,D,E,F,G,H}
{
  \fill (\P) circle (0.05);
}

\draw[thick, dashed, red] ($0.5*(B)+0.5*(C)$)--($0.5*(F)+0.5*(G)$);

\node at (90:1.4) {\scriptsize$E_{n+2}$};
\node at (135:1.4) {\scriptsize$E_1$};
\node at (175:1.4) {\scriptsize$\vdots $};
\node at (225:1.4) {\scriptsize$E_r $};
\node at (270:1.4) {\scriptsize$E_{r+1}$};
\node at (320:1.55) {\scriptsize$E_{r+2}$};
\node at (5:1.4) {\scriptsize$\vdots$};
\node at (40:1.55) {\scriptsize$E_{n+1}$};

\node at (2.5,0) {$\leadsto$};

\begin{scope}[shift={(5,0)}]
\coordinate (A) at (22.5:1.3);
\coordinate (B) at (67.5:1.3);
\coordinate (C) at (112.5:1.3);
\coordinate (D) at (157.5:1.3);
\coordinate (E) at (202.5:1.3);
\coordinate (F) at (247.5:1.3);
\coordinate (G) at (292.5:1.3);
\coordinate (H) at (337.5:1.3);
\coordinate (O) at (1/2,0);

\draw[thick, fill=yellow!30] (C)--(D)--(E)--(F)--(O)--cycle;
\foreach \P in {C,D,E,F,O}
{
  \fill (\P) circle (0.05);
}

\node at (135:1.5) {\scriptsize$E_1^+$};
\node at (175:1.4) {\scriptsize$\vdots $};
\node at (220:1.5) {\scriptsize$E_r^+$};
\node[right] at ($0.5*(F)+0.5*(O)$) {\scriptsize$E_{r+1}^+$};
\node[right] at ($0.5*(C)+0.5*(O)$) {\scriptsize$E_{r+2}^+$};
\node[right] at (O) {\scriptsize$v^+$};

\end{scope}

\begin{scope}[shift={(7.5,0)}]
\coordinate (A) at (22.5:1.3);
\coordinate (B) at (67.5:1.3);
\coordinate (C) at (112.5:1.3);
\coordinate (D) at (157.5:1.3);
\coordinate (E) at (202.5:1.3);
\coordinate (F) at (247.5:1.3);
\coordinate (G) at (292.5:1.3);
\coordinate (H) at (337.5:1.3);
\coordinate (O) at (-1/2,0);

\draw[thick, fill=yellow!30] (G)--(H)--(A)--(B)--(O)--cycle;
\foreach \P in {G,H,A,B,O}
{
  \fill (\P) circle (0.05);
}
\node at (320:1.55) {\scriptsize$E^-_{1}$};
\node at (5:1.4) {\scriptsize$\vdots$};
\node at (40:1.55) {\scriptsize$E^-_{s}$};
\node[left] at ($0.5*(B)+0.5*(O)$) {\scriptsize$E^-_{s+1}$};
\node[left] at ($0.5*(G)+0.5*(O)$) {\scriptsize$E^-_{s+2}$};

\node[left] at (O) {\scriptsize$v^-$};
\end{scope}

\end{tikzpicture}
\caption{Partition of a polygon.}
\label{fig_partition_polygon}
\end{figure}

If $\lambda(E_{n+2})$ and $\lambda(E_{r+1})$ are not parallel, then $\lambda^+$ and $\lambda^-$ are characteristic functions and we can construct their associated toric orbifolds $X(P^+,\lambda^+)$ and~$X(P^-,\lambda^-)$.
Writing $\pi_+$ and $\pi_-$ for their orbit maps, respectively, the toric orbifold $X(P, \lambda)$ is the $T^2$-equivariant connected sum 
\begin{equation}\label{eq_conn_sum}
X(P, \lambda) = X(P^+,\lambda^+)\# X(P^-,\lambda^-),
\end{equation}
obtained by removing small balls around the fixed points $\pi^{-1}_+(v^+)$ and $\pi^{-1}_-(v^-)$ and gluing along the resulting boundary sphere. We refer to  \cite[Construction 9.1.11]{BP-book} and \cite[\S 1.11]{DJ} for this construction. 

If $\lambda$ further satisfies that 
\begin{equation}\label{eq_conn_sum_det_pm1}
\lambda(E_{r+1})=(1,0)
\quad\text{and}\quad
\lambda(E_{n+2})=(0,1),
\end{equation}
then one can apply Theorem \ref{thm_main} to compute the cup products in $H^*(X(P,\lambda))$, as stated in the following theorem.

\begin{theorem}
Let $X(P,\lambda)$ be a $4$-dimensional toric orbifold where $P$ is an $(n+2)$-gon for $n\geq 2$ and $\lambda(E_i)=(a_i, b_i)$ for $i=1, \dots, n+2$. Suppose that there exists $r\in\{1,\ldots,n-1\}$ satisfying \eqref{eq_conn_sum_det_pm1}. Then, for $s=n-r$, there exists a basis 
\[
\{u_1,\ldots,u_r,w_1,\ldots,w_{s};v\}
\] 
for $\tilde{H}^*(X(P,\lambda))$
with $\deg u_i = \deg w_j =2$ and $\deg v=4$
such that their cup products are given by
\begin{align*}
\rm{(i)}&~u_i\cup u_j	=a_ib_jv &&\text{for } 1\leq i\leq j\leq r, \\ 
\rm{(ii)}&~ w_k\cup w_\ell	=-b_{r+1+k}a_{r+1+\ell}v && \text{for } 1\leq k\leq \ell\leq s, \\
\rm{(iii)}&~ u_p\cup w_q	=0 &&\text{for } 1\leq p\leq r \text{ and }1\leq q\leq s.
\end{align*}
\end{theorem}

\begin{proof}
Following the result of~\eqref{eq_conn_sum}, it suffices to compute the cup products on each of $H^*(X(P^+,\lambda^+))$ and $H^*(X(P^-,\lambda^-))$.

First, from Definition~\ref{dfn_connected sum} and Condition~\eqref{eq_conn_sum_det_pm1}, the characteristic function $\lambda^+$ is given by
\[
\lambda^+(E_i^+)=\begin{cases} 
\lambda(E_i) = (a_i, b_i) & \text{for } i=1, \dots, r;\\
\lambda(E_{r+1})=(1,0)  & \text{for } i=r+1;\\
\lambda(E_{n+2})=(0,1)  & \text{for } i=r+2.\\
\end{cases}
\]
Theorem~\ref{thm_main} implies that $H^*(X(P^+,\lambda^+))$ has a basis $\{u_1,\ldots,u_{r};v'\}$ such that 
\[
u_i\cup u_j=a_ib_jv'
\]
for $1\leq i\leq j\leq r$.

Next, the characteristic function $\lambda^-$ on $P^-$ is given by 
\[
\lambda^-(E_k^-)=\begin{cases} 
\lambda(E_{r+1+k}) = (a_{r+1+k}, b_{r+1+k}) & \text{for } k=1, \dots, s;\\
\lambda(E_{n+2})=(0,1)  & \text{for } k=s+1;\\
\lambda(E_{r+1})=(1,0)  & \text{for } k=s+2,\\
\end{cases}
\]
which does not directly satisfy the hypothesis of Theorem \ref{thm_main}. Hence, we consider a new characteristic function $\widetilde{\lambda}$ on $P^-$ by setting
\[
\widetilde{\lambda}(E_k^-)=\begin{cases} 
(b_{r+1+k}, a_{r+1+k}) & \text{for } k=1, \dots, s;\\
(1,0)  & \text{for } k=s+1;\\
(0,1)  & \text{for } k=s+2.\\
\end{cases}
\]
Then $\widetilde{\lambda}$ satisfies Equation~\eqref{eqn_lambda refined form} and Theorem~\ref{thm_main} implies that $H^*(X(P^-,\widetilde{\lambda}))$ has a basis $\{\tilde{w}_1,\ldots,\tilde{w}_{s};\tilde{v}\}$ such that
\[
\tilde{w}_k\cup\tilde{w}_l=b_{r+1+k}a_{r+1+\ell}\tilde{v}
\]
for $1\leq k\leq \ell\leq s$.
Consider the toric morphism
\[
X(id,\tau)\colon X(P^-,\lambda^-)\to X(P^-,\widetilde{\lambda})
\]
where $\tau\colon T^2\to T^2$ is given by $(x,y)\mapsto(y,x)$.
It is a homeomorphism, so it induces a basis $\{w_1,\ldots,w_{s};v''\}$ for $\tilde{H}^*(X(P^-,\lambda^-))$ where
\[
w_k=X(id,\tau)^*(\tilde{w}_k)
\quad\text{and}\quad
v''=X(id,\tau)^*(\tilde{v}).
\]
It follows from the naturality of cup products that
\[
w_k\cup w_l=b_{k+r+1}a_{\ell+r+1}v''
\]
for $1\leq k\leq \ell\leq s$.

Equip all toric orbifolds with their cellular orientations.
By Proposition~\ref{prop_fundamental class dual}, cohomology classes $v'$ and $\tilde{v}$ are the Kronecker duals of $[X(P^+,\lambda^+)]$ and $[X(P^-,\widetilde{\lambda})]$, while
\[
\langle v'',[X(P^-,\lambda^-)]\rangle=-1
\]
as $X(id,\tau)$ reverses orientation.
Let $v\in H^4(X(P,\lambda))$ be the generator dual to~$[X(P,\lambda)]$. Homeomorphism~\eqref{eq_conn_sum} then yields the cup products (i)--(iii).
\end{proof}

\end{document}